\documentclass[12pt,oneside,english]{amsart}
\usepackage[T1]{fontenc}
\usepackage[latin9]{inputenc}
\usepackage{comment}
\usepackage{babel}
\usepackage{verbatim}
\usepackage{mathrsfs}
\usepackage{amstext}
\usepackage{amsthm}
\usepackage{amssymb}
\usepackage{amsfonts}
\usepackage{amsmath}
\usepackage{esint}
\usepackage{enumerate}
\usepackage[unicode=true,pdfusetitle,
bookmarks=true,bookmarksnumbered=false,bookmarksopen=false,
breaklinks=false,pdfborder={0 0 1},colorlinks=false]
{hyperref}
\usepackage[margin=2.5cm]{geometry}
\usepackage[foot]{amsaddr}
\usepackage{mathtools}
\usepackage{ dsfont }
\usepackage{color,graphicx}
\usepackage[shortlabels]{enumitem}
\usepackage{array}
\usepackage[normalem]{ulem}

\makeatletter
\numberwithin{equation}{section}
\numberwithin{figure}{section}
\theoremstyle{plain}
\newtheorem{thm}{\protect\theoremname}[section]

\theoremstyle{definition}
\newtheorem{rem}[thm]{\protect\remarkname}
\theoremstyle{definition}
\newtheorem{defn}[thm]{\protect\definitionname}
\theoremstyle{plain}
\newtheorem{prop}[thm]{\protect\propositionname}
\theoremstyle{plain}
\newtheorem{lem}[thm]{\protect\lemmaname}
\theoremstyle{plain}

\theoremstyle{plain} 

\theoremstyle{definition}

\theoremstyle{definition}

\theoremstyle{definition}

\theoremstyle{definition}

\theoremstyle{definition}

\usepackage{tikz}
\usetikzlibrary{shapes,arrows}
\usepackage{verbatim}
\usepackage{amsthm}
\usepackage{amstext}

\DeclareMathOperator{\dist}{dist}

\DeclareMathOperator{\Leb}{Leb}

\DeclareMathOperator{\supp}{supp}

\DeclareMathOperator{\NN}{NN}
\DeclareMathOperator{\tol}{tol}

\newcommand{\R}{\mathbb R}
\newcommand{\Z}{\mathbb Z}
\newcommand{\N}{\mathbb N}

\newcommand{\PP}{\mathbb P}

\newcommand{\eps}{\varepsilon}

\newcommand{\mH}{\mathcal{H}}

\newcommand{\udim}{\overline{\dim}_B\,}

\newcommand{\hdim}{\dim_H}

\newcommand{\ud}{\overline{d}}
\newcommand{\ld}{\underline{d}}

\DeclareMathOperator{\Lip}{Lip}
\DeclareMathOperator{\Per}{Per}

\DeclareMathOperator{\rank}{rank}

\makeatother

\providecommand{\conjecturename}{Conjecture}
\providecommand{\corollaryname}{Corollary}
\providecommand{\definitionname}{Definition}
\providecommand{\examplename}{Example}
\providecommand{\lemmaname}{Lemma}
\providecommand{\problemname}{Problem}
\providecommand{\propositionname}{Proposition}
\providecommand{\remarkname}{Remark}
\providecommand{\theoremname}{Theorem}
\providecommand{\taskname}{Task}

\begin{document}

	\title{On the convergence of the False Nearest Neighbours algorithm}

	\author[A. \'{S}piewak]{Adam \'{S}piewak}
	\address{Institute of Mathematics, Polish Academy of Sciences, ul.~\'Sniadeckich 8, 00-656 Warszawa, Poland}
	\email{ad.spiewak@gmail.com}


	\begin{abstract}
		The False Nearest Neighbours (FNN) algorithm is one of the most widely used methods for choosing the embedding dimension, i.e.~the number $k$ of delay coordinates, in the time-delay reconstruction of a dynamical system from an observed scalar time series. We prove rigorous results on the behaviour of this algorithm for long time series. We consider a smooth dynamical system on a compact manifold $M$ and an invariant probability measure $\mu$. For a typical measurement function and $\mu$-almost every trajectory, we show that the proportion of false nearest neighbours tends to zero (provided that the tolerance threshold is large enough) once $k >\dim M$. On the other hand, for systems with super-polynomial decay of correlations, we prove that this proportion tends to one whenever $k$ is smaller than the Hausdorff dimension of $\mu$. Consequently, in this case, the FNN algorithm (run for long enough time series) returns an embedding dimension $k$ satisfying $\dim_H \mu \leq k \leq \dim M + 1$. This is below the threshold $2\dim M + 1$ required by the Takens embedding theorem, so the classical Takens embedding theory does not explain why reconstructions and predictions with $k$ chosen by the FNN algorithm work in practice. The explanation for this case is provided by a recent probabilistic theory of time-delay embeddings. We apply our results to the H\'enon and Lorenz attractors, showing that the FNN algorithm typically returns $k=2$ for the H\'enon map (in the chaotic case) and $k \leq 3$ for the Lorenz system.
	\end{abstract}

	\maketitle

	\section{Introduction}

	\subsection{Time-delay reconstruction and the choice of the embedding dimension}\label{subsec: intro takens}

	Time-delay reconstruction is one of the basic tools of nonlinear time series analysis \cite{PCFS80, T81, Abarbanel_book, KantzSchreiber03,BradleyKantz15, DatserisParlitz22}. Suppose that a deterministic dynamical system evolves according to $x_{i+1} = T(x_i)$, where $T : X \to X$ is a map on the phase space $X$, and that instead of the full state $x_i$ one can only record a single scalar quantity $h(x_i)$, where $h : X \to \R$ is a measurement function, also called an \emph{observable}. For a continuous-time system, $T$ is the map advancing the state by one sampling interval (or by the chosen time delay). Starting from an initial state $x = x_0$, one obtains a time series
	\begin{equation}\label{eq: time series h} h(x), \ldots, h(T^{n+k-1}(x)). \end{equation}
	The method of delays replaces the unknown states by \emph{time-delay vectors} consisting of $k$ consecutive measurements:
	\begin{equation}\label{def: z_i}
		y_i = y_i(x) = (h(T^i x), \ldots, h(T^{i+k-1}x)) \in \R^k, \qquad i = 0, \ldots, n.
	\end{equation}
	The number $k$ of delay coordinates is called the \textbf{delay length} or the \emph{embedding dimension}. If $k$ is large enough, one expects that the sequence $y_0, y_1, \ldots, y_n$ faithfully represents the original trajectory, so that the \emph{observed dynamics} $y_i \mapsto y_{i+1}$ in $\R^k$ can serve as a model of the system. Such reconstructions are routinely used for forecasting, noise reduction, estimation of dimensions and Lyapunov exponents, or detection of deterministic structure in data, see e.g.~\cite{FarmerSidorowich87, sm90nonlinear, hgls05distinguishing, SugiharaFishes} and the monographs \cite{Abarbanel_book, KantzSchreiber03}.

	The choice of the delay length is crucial. If $k$ is too small, the reconstruction folds the attractor onto itself: states which are far apart in the phase space may be represented by nearby (or even equal) time-delay vectors, and one cannot expect to reliably model the underlying dynamics or predict future values of \eqref{eq: time series h}. If $k$ is unnecessarily large, the analysis requires more data and becomes more sensitive to noise \cite{Abarbanel_book, KantzSchreiber03}. The classical Takens embedding theorem \cite{T81} is commonly invoked to give a mathematical answer to the question of how large $k$ should be: if the phase space is a compact manifold $M$, then for a typical system and a typical measurement function, the time-delay vectors determine the states uniquely as soon as $k > 2\dim M$. Sauer, Yorke and Casdagli \cite{SYC91} showed that for the reconstruction of a (possibly fractal) attractor (or of any compact invariant set) $X$, it is enough to take $k > 2\udim X$, where $\udim$ is the (upper) box-counting dimension. In practice, however, these dimensions are not known in advance, and the delay length has to be estimated from the time series itself. A number of methods for this task have been proposed, see e.g.~\cite{CenysPyragas88, LiebertPawelzikSchuster91, KBA92, Cao97, KennelAbarbanel02, Kraemer21}, and \cite{Abarbanel_book, KantzSchreiber03, KrakovkaFNN15} for an overview.

	\subsection{The False Nearest Neighbours algorithm}\label{subsec: FNN intro}

	One of the most commonly used methods is the classical \emph{False Nearest Neighbours (FNN) algorithm}, introduced by Kennel, Brown and Abarbanel \cite{KBA92} (see also \cite{Abarbanel_book}; similar ideas appeared earlier in \cite{CenysPyragas88, LiebertPawelzikSchuster91}). Its idea is as follows. For a candidate delay length $k$, one finds, for each time-delay vector $y_i$, its nearest neighbour $y_j$ among the other time-delay vectors in the sample. If $k$ is too small, then $y_j$ may be close to $y_i$ only because of the folding described above, and not because the corresponding states are close. Such \emph{false nearest neighbours} are detected by appending the next measurement to both time-delay vectors: if the difference $|h(T^{i+k}x) - h(T^{j+k}x)|$ between the next measurements is large compared to the distance $\|y_i - y_j\|$ (by a factor larger than a tolerance threshold $R_{\tol}$), or compared to the size of the attractor (by a factor larger than $A_{\tol}$), then $y_j$ is declared to be a false nearest neighbour of $y_i$. The algorithm starts with $k = 1$ and increases $k$ until the proportion of time-delay vectors with a false nearest neighbour drops to zero or becomes sufficiently small. The current value of $k$ is then returned as the final delay length. A precise description of the algorithm is given in Section \ref{subsec: FNN}. The FNN algorithm is implemented in standard software for nonlinear time series analysis \cite{TISEAN99} and is routinely used for choosing the delay length in applications, see e.g.~\cite{Oscilators, Streamflow07, Robots}. Its variants and refinements were studied in \cite{Cao97, HeggerKantz99, KennelAbarbanel02, NoisyFNN97, KrakovkaFNN15}, and the idea of false neighbours has also been used for detecting generalized synchronization \cite{CoupledSynchornization95}, for testing topological conjugacy of time series \cite{Dlotko_et_al} and for regularizing neural network reconstructions of attractors \cite{Gilpin20}. In this work we study rigorously the classical algorithm of  \cite{KBA92}. We expect that our methods can be useful also for studying its modifications.

	It has been observed for a number of dynamical systems with a known dimension of the phase space that the delay length returned by the FNN algorithm (or by similar methods) is considerably smaller than the bound $2\dim M + 1$ of the Takens theorem \cite{CenysPyragas88, LiebertPawelzikSchuster91, KBA92, Abarbanel_book}. For instance, for the $x$-coordinate of the Lorenz system, whose phase space is three-dimensional and whose attractor has dimension approximately $2.06$, the percentage of false nearest neighbours drops to zero already for $k = 3$ \cite[Section~4.2.1]{Abarbanel_book}, while the embedding theorems require at least $k=7$ (Takens \cite{T81}) or $k = 5$ (Sauer--Yorke--Casdagli \cite{SYC91}). As pointed out in \cite[Section~2.4.1]{Abarbanel_book}, the conditions of the embedding theorems are sufficient, but they need not be necessary in particular cases. Nevertheless, for delay lengths below these bounds the classical embedding theorems cannot be invoked to justify the validity of the reconstruction.

	The behaviour of the FNN algorithm has been studied extensively in numerical experiments, with particular attention paid to the choice of the thresholds, to the influence of noise and to the behaviour of the algorithm for stochastic data \cite{FredkinRice95, Cao97, NoisyFNN97, HeggerKantz99, KennelAbarbanel02, KrakovkaFNN15}. While the properties of time-delay embeddings have been intensively studied from a theoretical point of view (see e.g.~\cite{T81, SYC91, Stark99, SBDH03, EYWR18, TwistyTakens19, BGS20, BGSPredict, KoltaiKunde, BGOMY24, BGSLimits}), much less seems to be rigorously known about the behaviour of the algorithms which are being employed in applications. To the best of our knowledge, no rigorous mathematical results on the behaviour of the FNN algorithm have been obtained so far. In particular, no rigorous bounds on the delay length to which the algorithm converges have been given. The goal of this work is to provide such bounds.

	\subsection{Summary of the results}\label{subsec: summary}

	We give bounds on the delay length provided by the FNN algorithm for long time series generated by a typical measurement function along a typical trajectory of a smooth deterministic system on a compact manifold $M$. Here a \emph{typical measurement function} is understood in the sense of prevalence, see Section \ref{subsec: preval}. A \emph{typical trajectory} means the trajectory of $\mu$-almost every initial state, where $\mu$ is an invariant probability measure of the system. An important characteristic of the system turns out to be the Hausdorff dimension $\dim_H \mu$ of the measure $\mu$, see Section \ref{subsec: dim}. We prove that as the length of the time series tends to infinity, the following hold:
	\begin{enumerate}[(i)]
		\item If $k > \dim M$, then the proportion of time-delay vectors with a false nearest neighbour tends to zero, provided that the tolerance $R_{\tol}$ is large enough (Theorem \ref{thm: FNN diffeo}).
		\item If the system is sufficiently chaotic (in the sense of the super-polynomial decay of correlations) and $k < \dim_H \mu$, then the proportion of time-delay vectors with a false nearest neighbour tends to one, for any choice of the thresholds (Theorem \ref{thm: inverse FNN diffeo}, see also Theorem \ref{thm: inverse FNN general}).
	\end{enumerate}
	Consequently, for long time series generated by sufficiently chaotic systems, the FNN algorithm typically returns a delay length $k$ satisfying
	\begin{equation}\label{eq: FNN output range}
		\dim_H \mu \leq k \leq \dim M + 1.
	\end{equation}
	Whenever $k$ satisfies the inequalities above, the classical Takens embedding theorem cannot be invoked, as it requires $k > 2\dim M$, and may often be insufficient also for the Sauer--Yorke--Casdagli version \cite{SYC91} requiring $k > 2\udim X$. On the other hand,  one can apply in this case the recent probabilistic embedding theory \cite{BGS20, BGS22, BGSPredict, BGSLimits, SRegular25}, which builds on the conjectures of Schroer, Sauer, Ott and Yorke \cite{SSOY98} and requires only $k > \hdim \mu$. See Section \ref{subsec: prob Takens} for more details.

	Table \ref{table: regimes} summarizes these results, together with the known mathematical guarantees for the reconstruction in the respective ranges of the delay length.

	In certain cases, we can improve the upper bound in \eqref{eq: FNN output range} to $k \leq \dim M$. This holds if the measure $\mu$ generating the trajectory is concentrated on a compact invariant subset $X$ of zero volume on $M$ (see the precise statement in Theorem \ref{thm: FNN diffeo}). This allows for applications to strange attractors, and in fact our results can be applied to H\'enon and Lorenz attractors, providing a rigorous result supporting the numerical observations that the FNN algorithm typically returns $k=2$ for the H\'enon map (for the chaotic parameters studied by Benedicks and Carleson \cite{BC91}, for a typical point in its basin of attraction) and $k \leq 3$ for the Lorenz system. See Section \ref{sec: applications} for a detailed discussion.

	\begin{table}[ht]
		\renewcommand{\arraystretch}{1.35}
		\begin{tabular}{|l|l|>{\raggedright\arraybackslash}p{7cm}|}
			\hline
			\textbf{delay length} & \textbf{proportion of FNNs} & \textbf{guarantees for the reconstruction (typical observable)} \\ \hline
			$k < \dim_H \mu$ & {\renewcommand{\arraystretch}{1}\begin{tabular}[t]{@{}l@{}} $\to 1$ \ (Theorem \ref{thm: inverse FNN diffeo}) \\ \footnotesize(assuming super-polynomial \\ \footnotesize decay of correlations) \end{tabular}} & not one-to-one on any set of full measure \cite{BGSLimits} \\ \hline
			$\dim_H \mu < k \leq \dim M$ & {\renewcommand{\arraystretch}{1}\begin{tabular}[t]{@{}l@{}} not covered \\ \footnotesize(except for the case \\ \footnotesize $k = \dim M$ with \\ \footnotesize $\mH^k(\supp \mu) = 0$) \end{tabular}} &one-to-one on a set of full measure \cite{BGS20}, almost sure prediction \cite{BGSPredict} \\ \hline
			$\dim M < k \leq 2\dim M$ & $\to 0$ \ (Theorem \ref{thm: FNN diffeo}) & in addition: local diffeomorphism at almost every point \cite{SRegular25} \\ \hline
			$k > 2\dim M$ & $\to 0$ \ (Theorem \ref{thm: FNN diffeo}) & diffeomorphic embedding of the whole phase space \cite{T81, SYC91} \\ \hline
		\end{tabular}
		\vspace{0.2cm}
		\caption{Behaviour of the FNN algorithm for long time series generated by a typical observable along a typical trajectory, compared with the guarantees for the time-delay reconstruction (i.e.~the properties of the delay-coordinate map $\phi_{h,k}$).}
		\label{table: regimes}
	\end{table}

	\subsection{The False Nearest Neighbours algorithm: precise formulation}\label{subsec: FNN}

	We now describe the FNN algorithm precisely. We follow the exposition from \cite{Abarbanel_book}, using the notation of time-delay coordinate maps.  Given a map $T : X \to X$, an observable $h : X \to \R$ and $k \in \N$, let $\phi = \phi_h = \phi_{h,k} : X \to \R^k$ be the \textbf{$k$-delay coordinate map}
	\[ \phi_h(x) = (h(x), h(Tx), \ldots, h(T^{k-1} x)), \]
	so that the time-delay vectors \eqref{def: z_i} are given by $y_i = \phi_h(T^i x)$.

		\begin{defn}\label{defn: FNN}
			Fix parameters $R_{\tol} > 0$ and  $A_{\tol} > 0$. For given $x \in X,\ n \in \N$ and $0 \leq i <n$ let $0 \leq i_{\NN(x,n)} < n$ be such that\footnote{note that $y_i$ might have several nearest neighbours. As explained in Remark \ref{rem: tiebreaking}, the way of choosing amongst them does not influence our results.}
			\begin{equation}\label{eq: NN def excude} \|y_i - y_{i_{\NN(x,n)}}\| = \min\{ \|y_i - y_j\| : 0 \leq j < n,\ j \neq i \},
			\end{equation}
			where $y_i = \phi_h(T^i x)$	(so $y_{i_{\NN(x,n)}}$ is a \textbf{nearest neighbour} of $y_i$ among $(y_j)_{0 \leq j \neq i < n}$). We will also use the notation $\NN_{x,n}(\phi_h(T^i x)) = y_{i_{\NN(x,n)}}$. Further, let
			\[ R_A(x,n) := \left( \frac{1}{n + k} \sum \limits_{j=0}^{n + k-1} |h(T^j x) - \chi|^2 \right)^{1/2}, \text{ where } \chi = \frac{1}{n + k } \sum \limits_{j=0}^{n + k - 1} h(T^j x) \]
			and
			\[ \sigma(x,i,n) := \frac{|h(T^{i+k}x) - h(T^{i_{\NN(x,n)} + k}x)|}{R_A(x,n)}. \]
			We say that $y_i=\phi_h(T^i x)$ has a $(R_{\tol}, A_{\tol}, n)$\textbf{-false nearest neighbour} ($(R_{\tol}, A_{\tol}, n)$\textit{-FNN} for short) if\footnote{we interpret here $0/0 = 0$ and $x/0 = +\infty$ for $x > 0$. Therefore if $y_i = y_{i_{\NN}}$, then $\frac{|h(T^{i+k}x) - h(T^{i_{\NN} + k} x)|}{\|y_i - y_{i_{\NN}} \|} > R_{\tol}$ if and only if $h(T^{i+k}x) \neq h(T^{i_{\NN} + k} x)$. Similarly, if $R_A(x,n) = 0$, then $h(T^i x)$ is constant among $i=0, \ldots, n+k-1$, so $\sigma(x,i,n)  = 0$.}
			\begin{equation}\label{eq: FNN condition} \frac{|h(T^{i+k}x) - h(T^{i_{\NN(x,n)} + k} x)|}{\|y_i - y_{i_{\NN(x,n)}} \|} > R_{\tol} \text{ or } \sigma(x,i,n) > A_{\tol}. \end{equation}
		\end{defn}

		Let us now describe the classical \textit{false nearest neighbours algorithm} \cite{KBA92, Abarbanel_book}. The input is given by a single time series \eqref{eq: time series h}. Then one performs the following steps:
		\begin{enumerate}
			\item fix parameters $R_{\tol}$ and $A_{\tol}$,
			\item fix a delay length $k$, starting with $k=1$,
			\item\label{it: embed} embed the time series into $\R^k$ via the delay-coordinate map, i.e.~transform \eqref{eq: time series h} into \eqref{def: z_i},
			\item\label{it: check FNN} if there are no $(R_{\tol}, A_{\tol}, n)$-FNNs within the embedded time series \eqref{def: z_i}, then the algorithm outputs $k$ as the delay length,
			\item\label{it: increase k} otherwise replace $k$ with $k+1$ and repeat steps \eqref{it: embed}--\eqref{it: increase k}.
		\end{enumerate}
		
			In practice, the stopping criterion is often used in a relaxed form, i.e.~one increases $k$ until the percentage of false nearest neighbours drops to some acceptably small value, see e.g.~\cite{NoisyFNN97}. In this paper we study this relaxed version. More precisely, we fix an additional parameter $p \in (0,1)$ and replace \eqref{it: check FNN}  with
		
		\begin{enumerate}
			\item[(4')] if $\frac{1}{n} \# \left\{ 0 \leq i < n : y_i \text{ has a } (R_{\tol}, A_{\tol}, n)\text{-FNN}\right\} < p$, then the algorithm outputs $k$ as the delay length.
		\end{enumerate}

 The second condition in \eqref{eq: FNN condition}, involving $A_{\tol}$, accounts for nearest neighbours which are not close in the absolute sense, which is relevant for sparse data in high dimensions \cite[Section~4.2]{Abarbanel_book}.

	\subsection{Main results}\label{subsec: main results}
	We consider compact Riemannian $C^r$ manifolds with distance given by a Riemannian metric (one can take the Euclidean metric for manifolds embedded in $\R^N$; see \cite{LeeIntroSmoothManifolds} for a detailed introduction to smooth manifolds).
	
	Our first result deals with the problem of ensuring that the FNN algorithm will detect only a negligible proportion of false nearest neighbours if the delay length $k$ is large enough. Below, $\supp \mu$ denotes the topological support of a measure $\mu$, see Section \ref{sec: prelim}.

		\begin{thm}{\bf (Convergence of the False Nearest Neighbours algorithm)}\label{thm: FNN diffeo}
		Let $M$ be a compact $C^r, r \geq 1$ manifold and let $T : M \to M$ be a $C^r$-diffeomorphism. Let $\mu$ be a $T$-invariant Borel probability measure on $M$. If  one of the following holds:
		\begin{enumerate}[label=(A\arabic*), topsep=2pt]
			\item\label{it: A1} $k > \dim M$,
		\end{enumerate}
		\noindent or
		\begin{enumerate}[label=(A\arabic*), resume, topsep=2pt]
			\item\label{it: A2} $k \geq \dim M$ and $\mH^{\dim M}(\supp \mu) = 0$,
		\end{enumerate}
		\noindent then for a prevalent $C^r$ observable $h : M \to \R$, for  every $A_{\tol} > 0$
		\begin{equation}\label{eq: FNN to 0} \lim \limits_{R_{\tol} \to \infty }\ \mathop{\limsup}\limits_{n \to \infty} \frac{1}{n} \# \left\{ 0 \leq i < n : \phi_h(T^i x) \text{ has a } (R_{\tol}, A_{\tol}, n)\text{-FNN}\right\} = 0 \text{ for } \mu\text{-a.e. } x \in M.
		\end{equation}
	\end{thm}

		Here $\mH^{\dim M}$ denotes the $(\dim M)$-dimensional Hausdorff measure (see Definition \ref{def: hdim set}) with respect to the Riemannian distance on $M$, so the assumption $\mH^{\dim M}(\supp \mu) = 0$ means that $\supp \mu$ has zero volume. It holds in particular if $\hdim \supp \mu < \dim M$.

		Our second result shows that for systems with super-polynomial decay of correlations (see Definition \ref{defn: spdc}), if $k$ is small enough, then the FNN algorithm will detect an almost full proportion of false nearest neighbours.

			\begin{thm}{\bf (Divergence of the False Nearest Neighbours algorithm)}\label{thm: inverse FNN diffeo}
			Let $M$ be a compact $C^r, r \geq 1$ manifold and let $T : M \to M$ be a $C^r$-diffeomorphism. Let $\mu$ be a $T$-invariant Borel probability measure on $M$ having a super-polynomial decay of correlations. If $k < \hdim \mu$, then for a prevalent $C^r$ observable $h : M \to \R$, for every $R_{\tol} > 0,\ A_{\tol} > 0$
			\[ \lim \limits_{n \to \infty} \frac{1}{n} \# \left\{ 0 \leq i < n : \phi_h(T^i x) \text{ has a } (R_{\tol}, A_{\tol}, n)\text{-FNN}\right\} = 1 \text{ for } \mu\text{-a.e. } x \in M. \]
		\end{thm}

		Theorem \ref{thm: inverse FNN diffeo} is a special case of Theorem \ref{thm: inverse FNN general}, which holds for bi-Lipschitz maps on compact subsets of $\R^N$ and prevalent Lipschitz observables, see Section \ref{sec: general}.

		Combining Theorems \ref{thm: FNN diffeo} and \ref{thm: inverse FNN diffeo}, we see that for long time series generated by a typical observable along a typical orbit of a system as in Theorem \ref{thm: inverse FNN diffeo}, the FNN algorithm with the stopping criterion (4') for any $p \in (0,1)$ and a large threshold $R_{\tol}$ returns a delay length $k$ satisfying \eqref{eq: FNN output range}. Indeed, by Theorem \ref{thm: inverse FNN diffeo}, for every $k < \hdim \mu$ the proportion of false nearest neighbours tends to one, hence for $n$ large enough the algorithm does not stop at such $k$ (this holds also for the original stopping criterion \eqref{it: check FNN}), while by Theorem \ref{thm: FNN diffeo} applied with $k = \dim M + 1$, this proportion is smaller than any $p\in (0,1)$ for $n$ large enough, provided that $R_{\tol}$ is large enough.

		 \begin{rem}\label{rem: polynomial perturbations} In Theorems \ref{thm: FNN diffeo} and \ref{thm: inverse FNN diffeo}, prevalence is obtained concretely by random polynomial perturbations of the observable. Namely, embed $M$ into $\R^N$ (it is possible by the Whitney embedding theorem \cite{Whitney36}) and let $h_1, \ldots, h_m : \R^N \to \R$ be the basis of monomials of the linear space of real polynomials of $N$ variables of degree at most $D = 2k+1$. Then for every $C^r$ observable $h : M \to \R$, the conclusions of these theorems hold for the observable
		\[ h_\alpha = h + \sum_{j=1}^m \alpha_j h_j \]
		for Lebesgue-almost every $\alpha \in \R^m$, i.e.~for a random choice of the coefficients of the perturbing polynomial. Theorem \ref{thm: FNN diffeo} is true already with the choice $D=2k-1$ (see Section \ref{sec: FNN conv}).
		\end{rem}

		 {\begin{rem}\label{rem: tiebreaking} It follows from the proofs that, in the case of $y_i = \phi_h(T^i x)$ having several nearest neighbours, one can interpret the choice of the nearest neighbour to be identified as FNN in a way making both Theorems \ref{thm: FNN diffeo} and \ref{thm: inverse FNN diffeo} strongest. Namely, one can interpret ``\textit{$\phi_h(T^i x) \text{ has a } (R_{\tol}, A_{\tol}, n)\text{-FNN}$}'' as
			\begin{itemize}
			\item ``\textit{there exists a nearest neighbour of $\phi_h(T^i x)$ satisfying \eqref{eq: FNN condition}}'' in Theorem \ref{thm: FNN diffeo},
			\item ``\textit{all nearest neighbours of $\phi_h(T^i x)$ satisfy \eqref{eq: FNN condition}}'' in Theorem \ref{thm: inverse FNN diffeo}.
		\end{itemize}
		\end{rem}}

	\subsection[Relation to the time-delay embedding theory]{Relation to the time-delay embedding theory}\label{subsec: prob Takens}

	 Let us now discuss our results from the point of view of the time-delay embedding theory. The classical Takens embedding theorem requires $k > 2\dim M$ \cite{T81} (see also \cite{SYC91} for a version involving the dimension of the attractor). By Theorem \ref{thm: FNN diffeo}, for long time series the FNN algorithm typically returns $k \leq \dim M + 1 \leq 2\dim M$, i.e.~a delay length which is too small for the Takens theorem to be applicable. On the other hand, by Theorem \ref{thm: inverse FNN diffeo}, the returned delay length satisfies $k \geq \hdim \mu$, which is (up to the borderline case $k = \hdim \mu$) the condition required in the probabilistic embedding theory. More precisely, if $k > \hdim \mu$ and $\mu$ is ergodic, then for a typical observable $h$ the delay-coordinate map $\phi_h$ is injective on a set of full $\mu$-measure \cite[Theorem~1.2]{BGS20} (see also \cite[Theorem~4.3 and Remark~4.4]{BGS20}), and the prediction errors defined by Schroer, Sauer, Ott and Yorke \cite{SSOY98} converge to zero almost surely \cite[Theorem~1.5]{BGSPredict}. See \cite{BGSLimits} for a more detailed discussion of the results. This suggests that while the delay length $k$ obtained with the use of the FNN algorithm can be insufficient for a full reconstruction of the model (homeomorphically, without self-intersections) via the time-delay embedding, it is sufficient to guarantee that the reconstructed model has self-intersections only along a negligible (zero measure) set. Let us also highlight that from the algorithmic point of view, this yields some guarantees for the Farmer--Sidorowich prediction algorithm \cite{FarmerSidorowich87}. It aims at predicting the next time-delay vector from the previous ones, as the average of the successors of the past time-delay vectors close to the current one. For $x \in M$, $y \in \R^k$, $n \in \N$ and $\eps > 0$ let
	\[ \mathrm{Pred}_{x,n,\eps}(y) = \frac{1}{\# I_n} \sum_{i \in I_n} y_{i+1}, \quad \text{where } I_n = \{ 0 \leq i < n : \|y_i - y\| < \eps \} \]
	(provided $I_n \neq \emptyset$) and $y_i = \phi_h(T^i x)$ are the time-delay vectors \eqref{def: z_i}. If $\mu$ is ergodic and $k > \hdim \mu$, then for a prevalent observable $h$
	\[ \lim_{\eps \to 0} \lim_{n \to \infty} \mathrm{Pred}_{x,n,\eps}(\phi_h(y)) = \phi_h(Ty) \quad \text{for } \mu\text{-a.e. } x \in M \text{ and } \mu\text{-a.e. } y \in M, \]
	see \cite[Theorem~1.21 and Corollary~1.18]{BGSPredict}. In other words, whenever Theorem \ref{thm: inverse FNN diffeo} can be applied, the Farmer--Sidorowich algorithm with the delay length selected by the FNN algorithm correctly predicts the next time-delay vector at almost every point, provided that $k > \hdim \mu$ (which is automatically the case if $\hdim \mu$ is not an integer; otherwise it suffices to increase $k$ by one). This provides some theoretical explanation of the good performance of time-delay models with the delay length estimated by the FNN algorithm, observed in practice.

		\subsection{Discussion of the results and assumptions}\label{subsec: discussion}

		First, let us note that in Theorem \ref{thm: FNN diffeo} the threshold $R_{\tol}$ is taken large after passing with $n$ to infinity, while Theorem \ref{thm: inverse FNN diffeo} holds for all values of the thresholds. This is consistent with the observation that once enough data is available, the outcome of the algorithm is rather insensitive to the choice of the thresholds \cite[Section~4.2]{Abarbanel_book}, and that for an infinite amount of data the threshold $R_{\tol}$ can be chosen arbitrarily large \cite{NoisyFNN97}. Nevertheless, the problem of choosing $R_{\tol}$ and $A_{\tol}$ for finite samples is of practical importance and we do not address it.

		The assumption of super-polynomial decay of correlations is satisfied by a wide family of natural measures for chaotic dynamical systems, e.g.~by Gibbs measures on mixing basic sets of Axiom A diffeomorphisms, in particular by the SRB measures of mixing Axiom A attractors \cite{B75}. Exponential decay of correlations holds also for many non-uniformly hyperbolic systems \cite{Young98}, including H\'enon maps for a positive measure set of parameters \cite{BenedicksYoung00} and the Lorenz system \cite{AM16, AM17}.  It is however not clear to the author to what extent Theorem \ref{thm: inverse FNN diffeo} holds if the super-polynomial decay of correlations is not assumed. It is clear that some assumption is necessary: if $T = \mathrm{id}_M$ and $\mu$ is the normalized volume on $M$, then for every observable $h$ and any $x \in M$ all the time-delay vectors \eqref{def: z_i} coincide, so (with the conventions of Definition \ref{defn: FNN}) there are no false nearest neighbours and the algorithm returns $k = 1$, although $\hdim \mu = \dim M$. The proof of Theorem \ref{thm: inverse FNN diffeo} relies on recurrence properties of the original system and of the observed one (see Propositions \ref{prop: unpred} and \ref{prop: unrec}). Under the super-polynomial decay of correlations, the recurrence rate coincides almost surely with the Hausdorff dimension of $\mu$ \cite{Saussol06} (see also \cite{barreira2001hausdorff, RousseauSaussol10}). Therefore we expect that for general systems the critical threshold for the FNN algorithm is related  to both $\hdim \mu$ and the recurrence rate, but the author was unable to prove (or disprove) a more general result in this direction. The difficulty lies in the proof of the key Lemma \ref{lem: main ineq spdc}, which combines the decay of correlations, transversality bounds and the scaling coming from the dimension of the measure.

		Our results do not cover the intermediate range $\hdim \mu \leq k \leq \dim M$ (except for the case \ref{it: A2}), and they do not provide quantitative bounds on the length of the time series required for \eqref{eq: FNN output range} to hold. Even under the super-polynomial decay of correlations, it remains open to characterize $k$ in the range $[\hdim \mu, \dim M + 1]$ to which the FNN algorithm converges depending on the system. Finally, the behaviour of the algorithm for time series corrupted by noise (cf.~\cite{NoisyFNN97, HeggerKantz99}) is not studied.

		\section{Definitions and preliminaries}\label{sec: prelim}

		 A map $T : X \to X$ on a metric space $(X, d)$ is called \textbf{bi-Lipschitz} if there exists $L \geq 1$ such that $L^{-1} d(x,y) \leq d(Tx, Ty) \leq L d(x,y)$ for all $x, y \in X$. Note that such a map is injective, but it need not be surjective. For a Borel measure $\mu$ on $X$, we denote by $\supp \mu$ its topological support, i.e.~$\supp \mu = \{ x \in X : \mu(U) > 0 \text{ for every open neighbourhood $U$ of } x \}$.

		\subsection{Super-polynomial decay of correlations}

		For a Lipschitz function $\phi : X \to \R$ on a metric space $(X,d)$ we denote by $\Lip(\phi) = \sup \limits_{x \neq y} \frac{|\phi(x) - \phi(y)|}{d(x,y)}$ its Lipschitz constant and set $\|\phi\|_{\Lip} = \|\phi\|_\infty + \Lip(\phi)$.

		\begin{defn}\label{defn: spdc}
			Let $(X,T,\mu)$ be a probability measure preserving system on a metric space $X$. We say that $(X, T ,\mu)$ has a \textbf{super-polynomial decay of correlations} if for every Lipschitz $\phi, \psi : X \to \R$
			\[ \left| \int \phi \circ T^n \psi d\mu - \int \phi d\mu\int \psi d\mu \right| \leq \theta(n) \|\phi\|_{\Lip} \|\psi\|_{\Lip} , \]
			with $\theta$ satisfying $\lim \limits_{ n \to \infty} n^\alpha \theta(n) = 0$ for every $\alpha > 0$.
		\end{defn}
		
		\subsection{Dimension}\label{subsec: dim}

		{Let $X$ be a metric space.
		
		\begin{defn}\label{def: hdim set}
		For $s \geq 0$ and $A \subset X$, the \textbf{$s$-dimensional Hausdorff measure} of $A$ is
		\[ \mH^s(A) = \lim_{\delta \to 0} \inf \left\{ \sum_{i=1}^\infty (\mathrm{diam}\, U_i)^s : A \subset \bigcup_{i=1}^\infty U_i,\ \mathrm{diam}\, U_i \leq \delta \right\}, \]
		and the \textbf{Hausdorff dimension} of $A$ is $\hdim A = \inf\{ s \geq 0 : \mH^s(A) = 0 \}$.
		\end{defn}
		}

		For a measure $\mu$, we will use the following notions of dimensions.
		
		\begin{defn}
			Let $\mu$ be a finite Borel measure on $X$. The \textbf{lower} and \textbf{upper local dimensions} of $\mu$ at a point $x \in  \supp \mu$ are defined, respectively, as
			\[ \underline{d}(\mu, x) = \liminf_{r \to 0} \frac{\log \mu(B(x,r))}{\log r},\qquad  \overline{d}(\mu, x) = \limsup_{r \to 0} \frac{\log \mu(B(x,r))}{\log r}.\]
		\end{defn}
		
		\begin{defn}\label{defn: dim}
			The \textbf{(upper) Hausdorff dimension} of $\mu$ is defined as
			\[ \hdim \mu =  \underset{x \sim \mu}{\mathrm{ess\ sup}}\ \underline{d}(\mu, x). \]
		\end{defn}
		See \cite[Propositions 10.2--10.3]{FalconerTechniques} for characterizations in terms of Hausdorff dimensions of sets. We will use the following standard fact.
		
		\begin{lem}\label{lem: lower exact dim}
		Let $T : X \to X$ be a Lipschitz map and assume that $\mu$ has a super-polynomial decay of correlations. Then $\mu$ is ergodic and $\ld(\mu, x) = \hdim \mu$ for $\mu$-a.e. $x \in X$.
		\end{lem}
		\begin{proof} First, we prove ergodicity. Indeed, letting $n \to \infty$ in Definition \ref{defn: spdc} gives $\int \phi \circ T^n \psi d\mu \to \int \phi d\mu \int \psi d\mu$ for all Lipschitz $\phi, \psi : X \to \R$. As Lipschitz functions are dense in $L^2(\mu)$, the same convergence holds for all $\phi, \psi \in L^2(\mu)$, i.e.~$\mu$ is (strongly) mixing, hence ergodic (see e.g.~\cite[Theorem~1.23]{Walters82}). Therefore  $\ld(\mu, x)$ is constant almost surely by \cite[Proposition 10.6]{FalconerTechniques}, hence equal to its essential supremum $\hdim \mu$.
		\end{proof}

\subsection{Two auxiliary lemmas}\label{subsec: aux lemmas}

		Below, $\Leb$ denotes the Lebesgue measure on $\R^m$ restricted to the closed Euclidean ball $\overline{B_m}(0,1)$ and normalized so that it becomes a probability measure. For a linear map $\psi \colon \R^m \to \R^k$ (identified with its $k \times m$ matrix) we denote by $\sigma_1(\psi) \geq \sigma_2(\psi) \geq \ldots \geq \sigma_k(\psi) \geq 0$ its singular values, i.e.~the square roots of the eigenvalues (counted with multiplicity) of $\psi \psi^* \colon \R^k \to \R^k$, where $\psi^*$ is the adjoint of $\psi$ (see e.g.~\cite[Section~2.4]{MatrixComputations}).

		We will use the following lemma, proved as \cite[Lemma~4.2]{SYC91} (see also \cite[Lemma~14.3]{Rob11}).
		
		\begin{lem}\label{lem: key_ineq_inter}
			Let $\psi \colon  \R^m \to \R^k$ be a linear transformation. Assume $\sigma_p(\psi) > 0$ for some $p \in \{1, \ldots, k\}$. Then for every $z \in \R^{k}$ and $r > 0$,
			\[ \Leb(\{ \alpha \in B_m(0, 1) : \|\psi(\alpha) + z \| \leq r \}) \leq C \Big(\frac{r}{\sigma_p(\psi)}\Big)^p, \]
			where $C > 0$ depends only on $m,k$.
		\end{lem}
		
		We will also need the following modification.
		
		\begin{lem}\label{lem: eps delta geom bound}
		Let $\psi \colon  \R^m \to \R^{k+1}$ be a linear transformation such that $\sigma_{k+1}(\psi) > 0$. Then for every $z \in \R^{k+1}$ and $r, \delta > 0$ 
		\[ \Leb(\{ \alpha \in B_m(0, 1) : \psi(\alpha) + z \in [-r,r]^k \times [-\delta, \delta] \}) \leq C \frac{r^k \delta}{(\sigma_{k+1}(\psi))^{k+1}} , \]
		where $C > 0$ depends only on $m,k$.
		\end{lem}
	
		\begin{proof}
		{Let $Q = [-r,r]^k \times [-\delta, \delta]$ and $S = B_m(0,1) \cap \psi^{-1}(Q - z)$. As $\sigma_{k+1}(\psi) > 0$, $\psi$ has rank $k+1$, so $m \geq k+1$. By the coarea formula \cite[Lemma~3.4 and Theorem~3.6]{EvansGariepy},
		\[ \sqrt{\det(\psi \psi^*)} \, \mH^m(S) = \int_{\R^{k+1}} \mH^{m-k-1}(S \cap \psi^{-1}(\{y\})) \, d\mH^{k+1}(y), \]
		and $\sqrt{\det(\psi \psi^*)} = \prod_{i=1}^{k+1} \sigma_i(\psi) \geq \sigma_{k+1}(\psi)^{k+1} > 0$.
		For $y \notin Q - z$ the set $S \cap \psi^{-1}(y)$ is empty, while for $y \in Q - z$ it is contained in the intersection of $B_m(0,1)$ with an $(m-k-1)$-dimensional affine subspace, hence $\mH^{m-k-1}(S \cap \psi^{-1}(y)) \leq \mH^{m-k-1}(B_{m-k-1}(0,1))$. Consequently
		\[ \Leb(S) = \frac{\mH^m(S)}{\mH^m(B_m(0,1))} \leq \frac{\mH^{m-k-1}(B_{m-k-1}(0,1))}{\mH^m(B_m(0,1))} \cdot \frac{\mH^{k+1}(Q)}{\sqrt{\det(\psi \psi^*)}} \leq C \frac{r^k \delta}{\sigma_{k+1}(\psi)^{k+1}} \]
		with $C = 2^{k+1}\mH^{m-k-1}(B_{m-k-1}(0,1)) / \mH^m(B_m(0,1))$.
		}
		\end{proof}

		\subsection{Prevalence and the space of Lipschitz observables}\label{subsec: preval}
		Below we present the definition of prevalence -- a notion introduced by Hunt, Sauer and Yorke in \cite{Prevalence92} (an earlier version of this notion appeared in \cite{Christensen73}, see also \cite{BYbook00}), which may be considered as an analogue of the `Lebesgue almost sure' condition in infinite-dimensional normed linear spaces.
		
		\begin{defn}\label{defn:preval} By a (complete) \textbf{linear metric} in a linear space we mean a (complete) metric which makes the addition and scalar multiplication continuous. Let $V$ be a complete linear metric space (i.e.~a linear space with a complete linear metric). A Borel set $\mathcal S \subset V$ is called \textbf{prevalent} if there exists a Borel measure $\nu$ in $V$, which is positive and finite on some compact set in $V$, such that for every $v \in V$, we have $v + e \in \mathcal S$ for $\nu$-almost every $e \in V$. A non-Borel subset of $V$ is prevalent if it contains a prevalent Borel subset. For more information on prevalence we refer to \cite{Prevalence92} and \cite[Chapter~5]{Rob11}.
		\end{defn} 
		
		We will consider prevalence in the spaces $C^r(M, \R),\ r = 1, \ldots, \infty$ of $C^r$ observables as well as the space $\Lip(M, \R)$ of Lipschitz observables.

		{Prevalence in $C^r(M, \R)$ and $\Lip(X, \R)$ will be obtained by random polynomial perturbations of observables, similarly to \cite{SYC91, Rob11, BGS20, SRegular25}, in the following manner. Let $X \subset \R^N$ be a compact set, let $T : X \to X$ be a map and let $k \in \N$. Let $h_1, \ldots, h_m : \R^N \to \R$ be the basis of monomials of the linear space of real polynomials of $N$ variables of degree at most $D$ (we will set $D = 2k-1$ in Section \ref{sec: FNN conv} and $D = 2k+1$ in Section \ref{sec: inverse FNN}). As the measure $\nu$ in Definition~\ref{defn:preval} we employ $\nu = \xi\Leb$ for $\xi \colon \R^m \to \Lip(X,\R)$, $\xi(\alpha_1, \ldots, \alpha_m) = \sum_{j=1}^m \alpha_j h_j$, where $\Leb$ is the normalized Lebesgue measure on $\overline{B_m}(0,1)$. Then a set $\mathcal S \subset \Lip(X,\R)$ is prevalent if for every $h \in\Lip(X,\R)$, the function $h + \sum_{j=1}^m \alpha_j h_j$ is in $\mathcal S$ for $\Leb$-almost every $\alpha = (\alpha_1, \ldots, \alpha_m) \in \overline{B_m}(0,1)$. If $X$ is a manifold, then this also provides prevalence in the spaces $C^r(X, \R)$, $r=1,2,\ldots, \infty$. By \cite[Remark 2.3]{SRegular25}, it suffices to consider manifolds embedded in a Euclidean space.}

		{It will be useful to use the following notation. For a Lipschitz observable $h \colon X \to \R$ and $\alpha = (\alpha_1, \ldots, \alpha_m) \in \R^m$ let $h_\alpha \colon X \to \R$ be given by
		\[
		h_\alpha = h + \sum_{j=1}^m \alpha_j h_j.
		\]
		For simplicity, we write $\phi_\alpha$ instead of $\phi_{h_\alpha}$ for the $k$-delay coordinate map corresponding to $h_\alpha$, i.e.
		\[
		\phi_{\alpha} = \phi_{\alpha, k} \colon X \to \R^k, \qquad \phi_\alpha(x) = \phi_{\alpha, k}(x) = (h_\alpha(x), h_\alpha(Tx), \ldots, h_\alpha(T^{k-1}x)).
		\]
		{Note that for $x,y \in X$ we have
		\begin{equation}\label{eq:matrix_form} \phi_\alpha(x) - \phi_\alpha(y) = D_{x,y}\alpha + w_{x,y}
		\end{equation}
		for a $k\times m$ matrix $D_{x,y}$ defined by
		\begin{equation}\label{eq:D_xy}
			D_{x,y} = D_{x,y}(k) = \begin{bmatrix} h_1(x) - h_1(y) & \ldots & h_m(x) - h_m(y) \\
				h_1(Tx) - h_1(Ty) & \ldots & h_m(Tx) - h_m(Ty) \\
				\vdots & \ddots & \vdots \\
				h_1(T^{k-1}x) - h_1(T^{k-1}y) & \ldots & h_m(T^{k-1}x) - h_m(T^{k-1}y) \\
			\end{bmatrix}
		\end{equation}
		and
		\[
		w_{x,y} = \begin{bmatrix} h(x) - h(y) \\
			h(Tx) - h(Ty)\\
			\vdots \\
			h(T^{k-1} x) - h(T^{k-1}y) \end{bmatrix}.
		\]}}
		
		\subsection[Periodic points]{Periodic points}\label{subsec: per}

					For a map $T : X \to X$ and $x \in X$ denote
					
					\[ P_\ell(x) := \{ y \in X : \{y, Ty, \ldots, T^{\ell}y \} \cap \{x, Tx, \ldots, T^{\ell} x \} \neq \emptyset \}.  \]
					
					Moreover denote the set of periodic  points as
					\[ \Per_p(T) = \{ x \in X : T^p x = x \} \]
					for $p \in \N$.

We will need the following simple lemma, stated in the notation introduced above.
		\begin{lem}\label{lem: P disjoint k}
Assume that $T : X \to X$ is a bijection and fix  $x \in X$.  Then for $\Leb$-a.e.~$\alpha \in B_m(0,1)$
\[ \phi_\alpha(\{x\}) \cap \phi_\alpha(P_{k}(x) \setminus\{x\}) = \emptyset. \]
\end{lem}

\begin{proof}
{By \eqref{eq:matrix_form}, the set of $\alpha \in B_m(0,1)$ for which the conclusion fails is
\[ \bigcup \limits_{y \in P_k(x) \setminus \{x\}} \{ \alpha \in B_m(0,1) : \phi_\alpha(x) = \phi_\alpha(y) \} = \bigcup \limits_{y \in P_k(x) \setminus \{x\}} \{ \alpha \in B_m(0,1) : D_{x,y}\alpha = -w_{x,y} \}. \]
As $\rank D_{x,y} \geq 1$ for $y \neq x$ by \cite[Lemma 4.5]{BGS20}, each set in this union is contained in a proper affine subspace of $\R^m$, so it has zero Lebesgue measure. As $P_k(x)$ is finite, since $T$ is a bijection, this finishes the proof.}
\end{proof}

		\section{Convergence of the False Nearest Neighbours algorithm}\label{sec: FNN conv}

		In this section we prove Theorem \ref{thm: FNN diffeo}. The main tool is the following result on the regularity of time-delay embeddings with self-intersections, proved in \cite{SRegular25}. Given a compact $C^r$ manifold $M$, we fix a Riemannian metric on $M$ and denote by $\rho$ the induced distance on $M$ (every compact $C^r$ manifold admits a Riemannian metric, and since all such metrics are bi-Lipschitz equivalent, the choice of the metric does not matter for our purposes, see \cite[Remark~1.3]{SRegular25}). By the Whitney embedding theorem \cite{Whitney36}, we can assume that $M$ is a compact $C^r$ submanifold of some $\R^N$. As $M$ is compact, the Euclidean distance on $M$ is bi-Lipschitz equivalent to $\rho$, so the Hausdorff dimension and the condition $\mH^k(X) = 0$ do not depend on whether they are computed with respect to $\rho$ or to the Euclidean distance. We use the notation from Section \ref{subsec: preval}.

		\begin{thm}[{\cite[Theorem~4.1(1)]{SRegular25}}]\label{thm: SRegular}
			Let $M \subset \R^N$ be a compact $C^r$ submanifold, $r \geq 1$, and let $T : M \to M$ be a $C^r$ diffeomorphism. Let $X \subset M$ be a compact set and $\mu$ be a Borel probability measure on $X$. Fix $k \in \N$ and assume that $\mu \left( \bigcup \limits_{p=1}^{k-1} \Per_p(T) \right) = 0$. Let $h_1, \ldots, h_m : \R^N \to \R$ be a basis of the linear space of real polynomials of $N$ variables of degree at most $2k-1$. If $\mH^k(X) = 0$ and $k \geq \dim M$, then for every $C^r$ observable $h : M \to \R$ and Lebesgue-almost every $\alpha = (\alpha_1, \ldots, \alpha_m) \in \overline{B_m}(0,1)$, the following holds: for $\mu$-a.e. $x \in X$ there exists $C = C(x, \alpha) \in [0, \infty)$ so that
			\[ \|x - y\| \leq C \|\phi_{\alpha}(x) - \phi_{\alpha}(y)\| \text{ holds for every } y \in X. \]
		\end{thm}

	{We first prove Theorem \ref{thm: FNN diffeo} for ergodic measures, in the following form.

	\begin{thm}\label{thm: FNN ergodic}
	Let $M \subset \R^N$ be a compact $C^r$ submanifold, $r \geq 1$, let $T : M \to M$ be a $C^r$ diffeomorphism and let $\mu$ be an ergodic $T$-invariant Borel probability measure on $M$. Assume that $k \geq \dim M$ and $\mH^k(\supp \mu) = 0$, and let $h_1, \ldots, h_m$ be as in Theorem \ref{thm: SRegular}. Then for every $C^r$ observable $h : M \to \R$ and Lebesgue-almost every $\alpha \in \overline{B_m}(0,1)$ the following holds: for every $p > 0$ there exists $R = R(\alpha, p) > 0$ such that for every $R_{\tol} \geq R$ and every $A_{\tol} > 0$
	\[ \limsup \limits_{n \to \infty} \frac{1}{n} \# \left\{ 0 \leq i < n : \phi_\alpha(T^i x) \text{ has a } (R_{\tol}, A_{\tol}, n)\text{-FNN}\right\} < p \text{ for } \mu\text{-a.e. } x \in M. \]
	\end{thm}}

	\begin{proof}[{\bf Proof of Theorem \ref{thm: FNN ergodic}}]

	Let $X = \supp \mu$ and fix a $C^r$ observable $h : M \to \R$. If $\mu\left( \bigcup \limits_{\ell=1}^{k-1} \Per_{\ell}(T) \right) >0$, then by ergodicity $\mu$ is fully supported on a single $\ell$-periodic orbit for some $\ell \leq k-1$. By Lemma \ref{lem: P disjoint k}, for $\Leb$-a.e. $\alpha$, the map $\phi_\alpha$ is injective on this periodic orbit and the statement to be proved holds trivially in this case. We can therefore assume $\mu\left( \bigcup \limits_{\ell=1}^{k-1} \Per_{\ell}(T) \right) = 0$. By Theorem \ref{thm: SRegular}, for Lebesgue-almost every $\alpha \in \overline{B_m}(0,1)$ we have
	\[ \lim \limits_{L \to \infty} \mu(G_{L}) = 1, \]
	where
	\[ G_{L} = \{ z \in X : \underset{y \in X}{\forall} \ \|\phi_{\alpha}(Tz) - \phi_{\alpha}(Ty)\| \leq L\|\phi_{\alpha}(z) - \phi_{\alpha}(y)\|\}.\]
	Fix such $\alpha$. Fix $p>0$ and take $L = L(\alpha, p)$ such that $\mu(G_{L}) > 1 - p$. Fix $A_{\tol} > 0$. By the ergodic theorem, for $\mu$-a.e. $x \in X$
	\begin{equation}\label{eq: G_M freq} \lim \limits_{n \to \infty} \frac{1}{n} \# \left\{ 0 \leq i < n :  T^i x \in G_{L}\right\} > 1 - p.
	\end{equation}
	
		On the other hand, the ergodic theorem implies that
	\[ \lim \limits_{n \to \infty} R_A(x,n) = \sigma(h_{\alpha}):= \left( \int |h_{\alpha}(z) - \int h_{\alpha}(y)d\mu(y)|^2d\mu(z) \right)^{1/2} \text{ for } \mu\text{-a.e. } x \in X,\]
	so for $\mu$-a.e. $x \in X$, for all $n$ large enough
	\begin{equation}\label{eq: RA bound}
		R_A(x, n) \geq \sigma(h_{\alpha})/2.
	\end{equation}
	If $\sigma(h_{\alpha}) = 0$, then $h_{\alpha}$ is constant $\mu$-almost surely, and the conclusion of Theorem \ref{thm: FNN ergodic} holds for it trivially. Assume that $\sigma(h_{\alpha}) > 0$. As $\mu$ is ergodic, for $\mu$-a.e. $x$, every $\eps>0$ and every finite cover of $\supp \mu$ by $\eps$-balls, for all $n$ large enough, each element of the cover will be visited at least twice by the orbit $x, Tx, \ldots, T^{n-1}x$. As $\phi_{\alpha}$ is Lipschitz, this gives that for $\mu$-a.e. $x \in X$ and $n$ large enough, for every $0 \leq i < n$
	\begin{equation}\label{eq: NN dist} \|\phi_{\alpha}(T^i x) - \phi_{\alpha}(T^{i_{\NN(x,n)}}x) \| \leq A_{\tol}\sigma(h_{\alpha}) /2L.
	\end{equation}
	Take $R_{\tol} \geq L$. If $x \in X$ is such that $T^i x \in G_{L}$ and \eqref{eq: RA bound}, \eqref{eq: NN dist} hold, then (recall that $y_i = \phi_{\alpha}(T^i x)$)
	\[\begin{split} \frac{|h_{\alpha}(T^{i+k}x) - h_{\alpha}(T^{i_{\NN(x,n)} + k} x)|}{\|y_i - y_{i_{\NN(x,n)}} \|} &  \leq \frac{\|\phi_{\alpha}(T^{i+1} x) - \phi_{\alpha}(T^{i_{\NN(x,n)}+1}x) \|}{\|\phi_{\alpha}(T^i x) - \phi_{\alpha}(T^{i_{\NN(x,n)}}x) \|}  \leq L \leq R_{\tol}
	\end{split} \]
	and similarly
	\[\begin{split} \sigma(x,i,n) & = \frac{|h_{\alpha}(T^{i+k}x) - h_{\alpha}(T^{i_{\NN(x,n)} + k}x)|}{R_A(x,n)}  \leq \frac{ L\|\phi_{\alpha}(T^i x) - \phi_{\alpha}(T^{i_{\NN(x,n)}}x) \|}{R_A(x,n)} \\
		& \leq \frac{2L\|\phi_{\alpha}(T^i x) - \phi_{\alpha}(T^{i_{\NN(x,n)}}x) \|}{\sigma(h_{\alpha})} \leq A_{\tol}.
	\end{split}\]
	Therefore, we can conclude the following: for $\mu$-a.e. $x \in X$, for all $n$ large enough, for each $0 \leq i < n$
	\[ T^i x \in G_{L} \Longrightarrow \phi_{\alpha}(T^i x) \text{ does not have a } (R_{\tol}, A_{\tol}, n)\text{-FNN}.\]
	Recalling \eqref{eq: G_M freq} finishes the proof.

	\end{proof}

	{\begin{proof}[{\bf Proof of Theorem \ref{thm: FNN diffeo}}]
	As explained at the beginning of this section, we can assume that $M \subset \R^N$ is a compact $C^r$ submanifold. Let $h_1, \ldots, h_m$ be as in Theorem \ref{thm: SRegular} and let $\mu = \int \nu \, d\PP(\nu)$ be the ergodic decomposition of $\mu$ (see e.g.~\cite[Section~6.2]{Walters82}). For $\PP$-a.e.~$\nu$, the measure $\nu$ is an ergodic $T$-invariant Borel probability measure with $\nu(\supp \mu) = 1$, hence $\supp \nu \subset \supp \mu$. Consequently, $\mH^k(\supp \nu) = 0$ and $k \geq \dim M$: under \ref{it: A1} this holds as $\mH^k(M) = 0$ for $k > \dim M$, and under \ref{it: A2} as $\mH^{\dim M}(\supp \mu) = 0$ implies $\mH^k(\supp \mu) = 0$ for $k \geq \dim M$. Therefore Theorem \ref{thm: FNN ergodic} applies to $\PP$-a.e.~$\nu$. Whenever this is the case, \eqref{eq: FNN to 0} holds for $\Leb$-a.e. $\alpha \in B_m(0,1)$ and $\nu$-a.e. $x$. By Fubini's theorem, $\Leb$-a.e. $\alpha \in B_m(0,1)$ has the property that \eqref{eq: FNN to 0} holds for $\nu$-a.e. $x$ for $\PP$-a.e $\nu$, and hence for $\mu$-a.e $x$ as $\mu = \int \nu \, d\PP(\nu)$.
	\end{proof}}

	\section{Divergence of the False Nearest Neighbours algorithm}\label{sec: inverse FNN}

	\subsection[General version of Theorem \ref{thm: inverse FNN diffeo}]{General version of Theorem \ref{thm: inverse FNN diffeo}}\label{sec: general}

	In this section we prove the following general version of Theorem \ref{thm: inverse FNN diffeo}.

	\begin{thm}\label{thm: inverse FNN general}{\bf (Divergence of the False Nearest Neighbours algorithm)}\\
		Let $X \subset \R^N$ be a compact set and let $T : X \to X$ be a bi-Lipschitz map. Let $\mu$ be a $T$-invariant Borel probability measure on $X$ with super-polynomial decay of correlations. Let $k \in \N$ be such that $k < \hdim \mu$ and let $h_1, \ldots, h_m : \R^N \to \R$ be the basis of monomials of the linear space of real polynomials of $N$ variables of degree at most $2k+1$. Then for every Lipschitz observable $h : X \to \R$ and Lebesgue-almost every $\alpha = (\alpha_1, \ldots, \alpha_m) \in \overline{B_m}(0,1)$, for every $R_{\tol} >0,\ A_{\tol} > 0$
		\begin{equation}\label{eq: FNN one}
			\lim \limits_{n \to \infty} \frac{1}{n} \# \left\{ 0 \leq i < n : \phi_{\alpha}(T^i x) \text{ has a } (R_{\tol}, A_{\tol}, n)\text{-FNN}\right\} = 1 \text{ for } \mu\text{-a.e. } x \in X.
		\end{equation}
		Consequently, the conclusion holds for a prevalent Lipschitz observable $h : X \to \R$.
		\end{thm}

		Theorem \ref{thm: inverse FNN diffeo} is the special case of Theorem \ref{thm: inverse FNN general} for $C^r$-diffeomorphisms of compact $C^r$ manifolds (embedded in $\R^N$ by the Whitney embedding theorem).

		 Fix a Lipschitz observable $h : X \to \R$. The proof of Theorem \ref{thm: inverse FNN general} is based on the following two propositions. Below we use the notation $h_\alpha$, $\phi_\alpha$, $D_{x,y}$ and $w_{x,y}$ introduced in Section \ref{subsec: preval}.
					
					\begin{prop}\label{prop: unpred}
					Let $X \subset \R^N$ be a compact set, let $T : X \to X$ be a bi-Lipschitz bijection and let $\mu$ be a Borel probability measure on $X$ such that $\mu \left( \bigcup \limits_{p=1}^{k} \Per_p(T) \right) = 0$. Then for every $\beta \in (0,1]$, every $\eps, \eps'>0$ small enough (depending on $\beta$) and every $L > 0$, the following holds for $\Leb \otimes \mu$-a.e. $(\alpha, x) \in \overline{B_m}(0,1) \times X$: there exists $n_0 = n_0(\alpha, x, \beta, L, \eps, \eps')$ such that for every $j \in \Z$ with $|j| \geq n_0$ the following implication holds
					\begin{equation}\label{eq: unpred} \begin{rcases} \|\phi_\alpha(T^j x) - \phi_\alpha(x) \| \leq |j|^{-1/(k+\eps')} \\
					\|T^j x - x\| \geq |j|^{-\eps} \end{rcases} \Rightarrow \|\phi_\alpha(T^{j+1} x) - \phi_\alpha(Tx)\| \geq L |j|^{-\beta / (k+\eps')}
					\end{equation}
					\end{prop}
				
					\begin{prop}\label{prop: unrec}
					Under the assumptions of Theorem \ref{thm: inverse FNN general}, with $T$ being additionally a bijection, for every $\eps > 0$ and every $\eps'>0$ small enough (depending on $\eps$), the following holds for $\Leb \otimes \mu$-a.e. $(\alpha, x) \in \overline{B_m}(0,1) \times X$: there exists $n_1 = n_1(\alpha, x, \eps, \eps')$ such that for every $j \in \Z$ with $|j| \geq n_1$ the following implication holds
					\begin{equation}\label{eq: unrec}
					\| \phi_\alpha(T^j x) - \phi_\alpha(x)\| \leq |j|^{-1/(k+\eps')} \Rightarrow \|T^j x - x\| \geq |j|^{-\eps}.
					\end{equation}
					\end{prop}

Let us begin by explaining how to deduce Theorem \ref{thm: inverse FNN general} from Propositions \ref{prop: unpred} and \ref{prop: unrec}.
					
\begin{proof}[{\bf Proof of Theorem \ref{thm: inverse FNN general}}]
{We can assume that $T$ is a bijection, as otherwise we can restrict to the set $X_\infty = \bigcap_{n=0}^{\infty} T^n(X)$, which is compact, satisfies $T(X_\infty) = X_\infty$ and has full $\mu$-measure by the $T$-invariance of $\mu$.}
It suffices to fix $R_{\tol} > 1$ and $A_{\tol} > 0$ and show that \eqref{eq: FNN one} holds for a.e. $\alpha$ (as then one can take countable intersection over $R_{\tol}$ and $A_{\tol}$). Fix $L = 2R_{\tol}$ and $\beta = 1$ (in fact, the proof works for $\beta\in(0,1]$, see Remark \ref{rem: beta FNN}). Let $\eps, \eps'> 0$ be small enough so that the conclusions of both Propositions \ref{prop: unpred} and \ref{prop: unrec} are valid for them. Proposition \ref{prop: unpred} applies, as $\mu \left( \bigcup \limits_{p=1}^{\infty} \Per_p(T) \right) = 0$: otherwise the ergodicity of $\mu$ (see Lemma \ref{lem: lower exact dim}) would imply that $\mu$ is concentrated on a single periodic orbit, i.e.~on finitely many atoms, contradicting $k < \hdim \mu$.
					
Note that for given $x \in X$, applying \cite[Theorem~1.2]{BGS20} to the countable set $\{ T^i x : i \in \Z\}$ yields that for $\Leb$-a.e. $\alpha$
\begin{equation}\label{eq: orbit injective} \phi_\alpha \text{ is injective on the set } \{ T^i x : i \in \Z\}.
\end{equation}
By Fubini's theorem and Propositions \ref{prop: unpred}, \ref{prop: unrec}, $\Leb$-a.e. $\alpha \in \overline{B_m}(0,1)$ is such that for $\mu$-a.e. $x \in X$ property \eqref{eq: orbit injective} holds and \eqref{eq: unpred}, \eqref{eq: unrec} hold for $j \in \Z$ satisfying $|j| \geq n_0(x)= n_0(\alpha, x, \beta, L, \eps, \eps')$ and $|j| \geq n_1(x) = n_1(\alpha, x, \eps, \eps')$ (recall that $\beta, L, \eps, \eps'$ are already fixed). Let us fix such $\alpha$.
					
By \cite[Theorem 2]{RousseauSaussol10} (applied to the systems $(X,T,\mu)$ and $(X, T^{-1}, \mu)$) and the inequality $\ud(\phi_\alpha \mu, \phi_\alpha x) \leq k$ valid for $\mu$-a.e. $x \in X$ (see \cite[Chapter~10.1]{FalconerTechniques}), we have that for $\mu$-a.e. $x \in X$ there exists $n_2(x) = n_2(\eps', \alpha, x)$ such that
\begin{equation}\label{eq: rec n2}
\dist(\{\phi_\alpha(x) \}, \{ \phi_\alpha(Tx), \phi_\alpha(T^2 x), \ldots, \phi_\alpha(T^{\lfloor n/2 \rfloor}x)\}) \leq n^{-1/(k+\eps')} \text{ for } n \geq n_2(x)
\end{equation}
and
\begin{equation}\label{eq: rec n2 inv}
\dist(\{\phi_\alpha(x) \}, \{ \phi_\alpha(T^{-1}x), \phi_\alpha(T^{-2} x), \ldots, \phi_\alpha(T^{ - \lfloor n/2 \rfloor}x)\}) \leq n^{-1/(k+\eps')} \text{ for } n \geq n_2(x).
\end{equation}
Fix an arbitrary $\eps_0 > 0$ and let $M \in \N$ be such that
\[ \mu\left( F \right) \geq 1 - \eps_0, \text{ where } F = \left\{ x \in X : n_0(x), n_1(x), n_2(x) \leq M\right\}.\]
For $x \in F$ let
\[ \kappa(x) = \min \left\{ \| \phi_\alpha(x) - \phi_\alpha(T^j x)\| : j \in \Z \setminus \{ 0\}, |j| \leq M\right\}.\]
As $\mu(\Per_p(T)) = 0$ for every $p \in \N$, we see by \eqref{eq: orbit injective} that $\kappa(x) > 0$ for $\mu$-a.e. $x \in F$. Fix $\kappa_0 > 0$ such that
\begin{equation}\label{eq: F'} \mu(F') \geq 1- 2\eps_0, \text{ where } F' = \{  x \in F : \kappa(x)\geq \kappa_0\}.
\end{equation}
Let $M_0 \geq M$ be such that
\begin{equation}\label{eq: M0 choice} M_0^{-1/(k + \eps')} < \kappa_0/2 \text{ and } M_0^{-(1 - \beta)/(k + \eps')} < R_{\tol}
\end{equation}
(note that if $\beta = 1$, then the second condition is fulfilled trivially as we assume $R_{\tol} > 1$). Observe the following: for every $x \in X$,  every $n \geq M_0$ and every $0 \leq i < n$,
\begin{equation}\label{eq: F' gives distance}  \text{ if } T^i x \in F', \text{ then}\footnote{this holds for \emph{any} choice of a nearest neighbour $\phi_\alpha(T^{i_{\NN(x,n)}} x)$ of $\phi_\alpha(T^i x)$}\ |i_{\NN(x,n) }- i| > M \geq n_0(T^i x), n_1(T^i x), n_2(T^i x).
\end{equation}
Indeed, if $T^i x \in F'$ and $n \geq M_0$, then by \eqref{eq: rec n2} and \eqref{eq: rec n2 inv} we have
\begin{equation}\label{eq: FNN rec}\|\NN_{x,n}(\phi_\alpha(T^i x)) - \phi_\alpha(T^i x)\| \leq n^{-1 / (k + \eps')}.
\end{equation}
To see this, note that as $T^i x \in F' \subset F$, we have $n \geq M_0 \geq M \geq n_2(T^i x)$, so \eqref{eq: rec n2} and \eqref{eq: rec n2 inv} can be applied with $T^i x$ in place of $x$. If $i + \lfloor n/2 \rfloor < n$, then \eqref{eq: rec n2} gives $j \in \{ i+1, \ldots, i + \lfloor n/2 \rfloor \} \subset \{0, \ldots, n-1\}$ such that $\|\phi_\alpha(T^j x) - \phi_\alpha(T^i x)\| \leq n^{-1/(k+\eps')}$. Otherwise $i \geq n - \lfloor n/2 \rfloor \geq \lfloor n/2 \rfloor$, and \eqref{eq: rec n2 inv} gives such $j \in \{ i - \lfloor n/2 \rfloor, \ldots, i-1 \} \subset \{0, \ldots, n-1\}$. In both cases $j \neq i$, so \eqref{eq: FNN rec} follows from the definition of the nearest neighbour \eqref{eq: NN def excude}.
 As $n \geq M_0$, this gives
\[ \|\NN_{x,n}(\phi_\alpha(T^i x)) - \phi_\alpha(T^i x)\| \leq M_0^{- 1 / (k +\eps')} < \kappa_0/2 < \kappa(T^i x)\]
and consequently \eqref{eq: F' gives distance} holds by the definition of $\kappa(T^i x)$.
{\par\medskip\noindent\textbf{Claim.} For every} $n \geq M_0$ and  $0 \leq i < n$
\begin{equation}\label{eq: F' gives FNN} \text{ if } T^i x \in F', \text{ then } \phi_\alpha(T^i x) \text{ has a } (R_{\tol}, A_{\tol}, n)\text{-FNN}.
\end{equation}
{\noindent\textbf{Proof of the claim.}} \eqref{eq: F' gives distance} implies that if $T^i x \in F'$, then any nearest neighbour $\NN_{x,n}(\phi_\alpha(T^i x)) = \phi_\alpha(T^{i_{\NN(x,n) }}x)$ satisfies $|i_{\NN(x,n) }- i| \geq n_1(T^i x)$, while \eqref{eq: FNN rec} implies
\begin{equation}\label{eq: FFN rec i}\|\NN_{x,n}(\phi_\alpha(T^i x)) - \phi_\alpha(T^i x)\| \leq |i_{\NN(x,n) }- i|^{- 1 / (k + \eps')}.
\end{equation}
Therefore we can invoke \eqref{eq: unrec} to claim that
\begin{equation}\label{eq: FFN rec X}\|T^{i_{\NN(x,n)}} x - T^i x\| \geq |i_{\NN(x,n) }- i|^{-\eps}.
\end{equation}
As \eqref{eq: F' gives distance} implies also $|i_{\NN(x,n) }- i| \geq n_0(T^i x)$ for $T^i x \in F'$,  \eqref{eq: unpred} together with \eqref{eq: FFN rec i}, \eqref{eq: FFN rec X} give
\[   \|\phi_\alpha(T^{i_{\NN(x,n)} + 1} x) - \phi_\alpha(T^{i+1} x) \| \geq L |i_{\NN(x,n) }- i|^{-\beta / (k+\eps')}. \]
Combined with \eqref{eq: FFN rec i} this gives
\begin{equation}\label{eq: FNN delay 2R} \frac{\|\phi_\alpha(T^{i_{\NN(x,n)} + 1} x) - \phi_\alpha(T^{i+1} x)\|}{\|\phi_\alpha(T^i x) - \NN_{x,n}(\phi_\alpha(T^i x))\|^\beta } \geq L = 2R_{\tol}.
\end{equation}
Moreover, by \eqref{eq: M0 choice} and \eqref{eq: FNN rec}
\[ \|\NN_{x,n}(\phi_\alpha(T^i x)) - \phi_\alpha(T^i x)\|^{1 - \beta} < R_{\tol}.\]
{Further, the first $k-1$ coordinates of $\phi_\alpha(T^{i_{\NN(x,n)} + 1} x) - \phi_\alpha(T^{i+1} x)$ are the same as the last $k-1$ coordinates of $\phi_\alpha(T^{i_{\NN(x,n)}} x) - \phi_\alpha(T^i x)$, while its last coordinate is $h_\alpha(T^{i_{\NN(x,n)} + k} x) - h_\alpha(T^{i+k} x)$, so
\begin{equation}\label{eq: last coord}
\begin{split}
\|\phi_\alpha(T^{i_{\NN(x,n)} + 1} x) - \phi_\alpha(T^{i+1} x)\| \leq\ & |h_\alpha(T^{i_{\NN(x,n)} + k} x) - h_\alpha(T^{i+k} x)| \\
& + \|\phi_\alpha(T^{i_{\NN(x,n)}} x) - \phi_\alpha(T^i x)\|.
\end{split}
\end{equation}}
Hence \eqref{eq: FNN delay 2R} and \eqref{eq: last coord} yield, as $\beta = 1$ and $R_{\tol} > 1$, that $\frac{|h_\alpha(T^{i+k} x) - h_\alpha(T^{i_{\NN(x,n)} + k} x)|}{\|y_i - y_{i_{\NN(x,n)}}\|} \geq 2R_{\tol} - 1 > R_{\tol}$ for $y_i = \phi_\alpha(T^i x)$, so  that indeed $\phi_\alpha(T^i x)$ has a $ (R_{\tol}, A_{\tol}, n)$-FNN, proving \eqref{eq: F' gives FNN} (note that the criterion involving $A_{\tol}$ is not used here at all). This finishes the proof of the claim.
{\par\medskip}
By the ergodic theorem (recall that $\mu$ is ergodic by Lemma \ref{lem: lower exact dim})
\[ \liminf \limits_{n \to \infty} \frac{1}{n} \sum \limits_{i=0}^{n-1} \mathds{1}_{F'}(T^i x) \geq \mu(F') \text{ for } \mu\text{-a.e. } x \in X,\]
so \eqref{eq: F'} and \eqref{eq: F' gives FNN} give
\[ \liminf \limits_{n \to \infty} \frac{1}{n} \# \left\{ 0 \leq i < n : \phi_\alpha(T^i x) \text{ has a } (R_{\tol}, A_{\tol}, n)\text{-FNN}\right\} \geq 1-2\eps_0 \text{ for } \mu\text{-a.e. } x \in X. \]
As $\eps_0$ was arbitrary, the proof is complete.
\end{proof}

{\begin{rem}\label{rem: beta FNN}
It is easy to see that $\beta \in (0,1]$ can be chosen arbitrarily in the above proof. This gives that under the assumptions of Theorem \ref{thm: inverse FNN general}, the same conclusion holds for the version of the FNN algorithm in which condition \eqref{eq: FNN condition} in Definition \ref{defn: FNN} is replaced by the following one: $y_i$ is said to have a $(\beta, R_{\tol}, A_{\tol}, n)$-FNN if
\[ \frac{|h(T^{i+k}x) - h(T^{i_{\NN(x,n)} + k} x)|}{\|y_i - y_{i_{\NN(x,n)}} \|^\beta} > R_{\tol} \text{ or } \sigma(x,i,n) > A_{\tol}. \]
It corresponds to requiring the next step prediction to be only $\beta$-H\"older rather than Lipschitz as in \eqref{eq: FNN condition}.
\end{rem}}

Let us now turn to the proofs of Propositions \ref{prop: unpred} and \ref{prop: unrec}.

\subsection[Proof of Proposition \ref{prop: unpred}]{Proof of Proposition \ref{prop: unpred}}\label{subsec: proof unpred}

First, we will need Lemma \ref{lem: sigma bound k} below and Lemma \ref{lem: P disjoint k}, the first of which concerns also the matrix $D_{x,y}(k+1)$ (this is why we require $h_1, \ldots, h_m$ to be the basis of monomials of the space of polynomials of degree at most $2k+1$ rather than $2k-1$).

\begin{lem}\label{lem: sigma bound k} Let $X \subset \R^N$ be a compact set and let $T : X \to X$ be a bi-Lipschitz bijection. Fix  $x \in X \setminus \bigcup \limits_{p=1}^{k} \Per_p(T)$. There exists $\xi = \xi(x) > 0$ with the following property. For every $\eta > 0$ there exists $C = C(x, \eta)$ such that
	\begin{equation}\label{eq: sigma bound close}
\sigma_k(D_{x,y}(k)) \geq \sigma_{k+1}(D_{x,y}(k+1)) \geq C\|x - y\| \ \text{ for } y \in X \text{ with } \|x - y\| \leq \xi,
\end{equation}
and
\begin{equation}\label{eq: sigma bound far}
\sigma_{k+1}(D_{x,y}(k+1)) \geq C \ \text{ for } y \in X \text{ with } \dist(y, P_{k}(x)) \geq \eta.
\end{equation}
\end{lem}

\begin{proof}
The first inequality in \eqref{eq: sigma bound close} holds because $D_{x,y}(k)$ consists of the first $k$ rows of $D_{x,y}(k+1)$, and $\sigma_{k+1}(D_{x,y}(k+1)) = \min \{ \|D_{x,y}(k+1)^T u\| : u \in \R^{k+1},\ \|u\| = 1 \}$, while restricting this minimum to the vectors $u$ with the last coordinate equal to zero gives $\sigma_k(D_{x,y}(k))$. In the rest of this proof we write $D_{x,y} = D_{x,y}(k+1)$. As $x \notin \bigcup \limits_{p=1}^{k} \Per_p(T)$ and $T$ is a Lipschitz bijection, there exist $\xi = \xi(x) > 0$ and $\eps_0 = \eps_0(x)$ such that for every $y \in X$ satisfying $\|x - y\| \leq \xi$, inequalities
\[ \|T^i z_1 - T^j z_2\| \geq \eps_0 \text{ for } i,j=0, \ldots, k \text{ such that } i \neq j,\ z_1, z_2 \in \{x,y\} \]
hold. Therefore \cite[Proposition 4.4]{BGSLimits} implies that there exists a constant $C$ such that
\[ \sigma_{k+1}(D_{x,y}) \geq C\eps_0^{2k}\|T^{k} x - T^{k}y\| \text{ if } \|x - y\| \leq \xi.\]
As $T$ is a bi-Lipschitz bijection, $T^{-k}$ is Lipschitz, which gives the second inequality in \eqref{eq: sigma bound close}, after decreasing $C$ if necessary. We will finish the proof by showing that
\begin{equation}\label{eq:rank Dxy} \rank D_{x,y} = k+1 \text{ if } y \notin P_{k}(x).
\end{equation}
This will imply $\sigma_{k+1}(D_{x,y}) > 0$ on a compact set $\{ y \in X : \dist(y, P_{k}(x)) \geq \eta \}$, so  $\sigma_{k+1}(D_{x,y}) \geq C = C(x,\eta) > 0$ on this set, by the continuity of the $(k+1)$-th singular value  (see e.g.~\cite[Corollary 8.6.2]{MatrixComputations}).
To prove \eqref{eq:rank Dxy}, denote first the elements of the set $\{x, Tx, \ldots, T^{k}x, y, Ty, \ldots, T^{k}y \}$ by $z_1, \ldots, z_{\ell}$ (without multiplicities), preserving the order within the set, so that
\[\ell = \#\{x, Tx, \ldots, T^{k}x, y, Ty, \ldots, T^{k}y \}.\]
Under this notation, the matrix $D_{x,y}$ can be written as a product of matrices
\[ D_{x,y} = J_{x,y} V_{x,y}, \]
where 
\[ V_{x,y}
= \begin{bmatrix} h_1(z_1) & \ldots & h_m(z_1) \\
	\vdots & \ddots & \vdots \\
	h_1(z_{\ell}) & \ldots & h_m(z_{\ell})
\end{bmatrix} \]
and $J_{x,y} = J_{x,y}(T,k+1) = [m_{i,j}]_{\substack{1 \leq i \leq k+1\\1 \leq j \leq \ell}}$ is a $(k+1) \times \ell$ matrix, where
\[m_{i,j} = \begin{cases} 1 &\text{if}\ z_j = T^{i-1} x\\
	-1 &\text{if}\ z_j = T^{i-1} y\\
	0 &\text{otherwise}
\end{cases}.\]
Since $z_1, \ldots, z_{\ell}$ are distinct, $\ell \leq 2(k+1)$ and polynomials of $N$ variables of degree at most $2k+1$ can attain arbitrary values at any $2(k+1)$ distinct points of $\R^N$ (see e.g.~\cite[Section~1.2, eq.~(1.9)]{poly-interpolation}), the matrix $V_{x,y}$ has full rank. Therefore, $\rank D_{x,y} = \rank J_{x,y}$. As $T$ is a bijection, if $x \notin \bigcup \limits_{p=1}^{k} \Per_p(T)$ and $y \notin P_{k}(x)$, then $\{z_1, \ldots, z_{k+1}\} = \{x, Tx, \ldots, T^{k}x\}$ and $\{y, Ty, \ldots, T^{k}y\} \cap \{z_1, \ldots, z_{k+1}\} = \emptyset$, hence we see that $J_{x,y}$ is of the block form
\[ J_{x,y} = \begin{bmatrix} I & A
\end{bmatrix}, \]
where $I$ is a $(k+1) \times (k+1)$ identity matrix and $A$ is a $(k+1) \times (\ell - k - 1)$ matrix. This shows $\rank J_{x,y} = k+1$ in this case and finishes the proof.
\end{proof}

\medskip
For $\eta > 0$ and $x \in X \setminus \bigcup \limits_{p=1}^{k} \Per_p(T)$ with $\xi(x)$ as in Lemma \ref{lem: sigma bound k}, define
\[ K_\eta(x) = \{ y \in X : \| x - y\| \leq \xi(x) \text{ or } \dist(y, P_{k}(x)) \geq \eta \}. \]

\begin{lem}\label{lem: k+1 leb bound}
Fix $x \in X \setminus \bigcup \limits_{p=1}^{k} \Per_p(T)$. Then for every $\eta > 0$ there exists $C = C(x,\eta)$ such that for every $y \in K_\eta(x) \setminus \{x\}$ inequality
\[\Leb \left( \left\{ \alpha \in B_m(0,1) : \|\phi_\alpha(x) - \phi_\alpha(y)\| \leq r, \|\phi_\alpha(Tx) - \phi_\alpha(Ty)\| \leq \delta \right\} \right) \leq C \frac{r^k \delta}{\| x - y\|^{k+1}}\]
holds for every  $r, \delta > 0$.
\end{lem}

\begin{proof}
By Lemma \ref{lem: eps delta geom bound}
\[\begin{split}
\Leb & \left( \left\{ \alpha \in B_m(0,1) : \|\phi_\alpha(x) - \phi_\alpha(y)\| \leq r, \|\phi_\alpha(Tx) - \phi_\alpha(Ty)\| \leq \delta \right\} \right) \\
& \leq C \frac{r^k \delta}{\sigma_{k+1}(D_{x,y}(k+1))^{k+1}}.
\end{split}
\]
Combining this with \eqref{eq: sigma bound close} and \eqref{eq: sigma bound far} and using the fact that $X$ is compact  finishes the proof.

\end{proof}

\begin{lem}\label{lem: keta unpred}
For every $\beta \in (0,1]$, every $\eps, \eps'>0$ small enough (depending on $\beta$), the following holds. For every $x \in X \setminus \bigcup \limits_{p=1}^{k} \Per_p(T)$, for  $\Leb$-a.e. $\alpha \in \overline{B_m}(0,1)$, for every $\eta, L > 0$, there exists $n_0 = n_0(\alpha, x, \beta, L, \eps, \eps', \eta)$ such that for every $j \in \Z$ with $|j| \geq n_0$ the following implication holds
\begin{equation*} \begin{rcases} \|\phi_\alpha(T^j x) - \phi_\alpha(x) \| \leq |j|^{-1/(k+\eps')} \\
		\|T^j x - x\| \geq |j|^{-\eps} \\
	  T^j x \in K_\eta(x) \end{rcases} \Rightarrow \|\phi_\alpha(T^{j+1} x) - \phi_\alpha(Tx)\| \geq L |j|^{-\beta / (k+\eps')}.
\end{equation*}
\end{lem}

\begin{proof}
Fix $\eta, L > 0$. If  $T^j x \in K_\eta(x)$ and $\|T^j x - x\| \geq |j|^{-\eps}$, then by Lemma \ref{lem: k+1 leb bound}
\[\begin{split}
\Leb & \Big( \Big\{ \alpha \in B_m(0,1) : \|\phi_\alpha(T^j x) - \phi_\alpha(x)\| \leq |j|^{- 1 / (k + \eps')}, \\
& \qquad \|\phi_\alpha(T^{j+1}x) - \phi_\alpha(Tx)\| \leq L |j|^{-\beta / (k+\eps')} \Big\} \Big) \\
& \leq L C(x,\eta) |j|^{- \frac{k + \beta}{k + \eps'}}\|T^j x - x\|^{-(k+1)} \\
& \leq L C(x,\eta) |j|^{- \frac{k + \beta}{k + \eps'} + \eps(k+1)}.
\end{split} \]
For given $\beta>0$, if we choose $\eps, \eps'>0$ small enough, then $- \frac{k + \beta}{k + \eps'} + \eps(k+1) < -1$, hence summing the above probabilities over $\left\{ j \in \Z : \|T^j x - x\| \geq |j|^{-\eps},\ T^j x \in K_\eta(x)\right\}$ and applying the Borel--Cantelli lemma gives the desired property for fixed $\eta, L >0$. Taking countable intersection over $\eta \to 0$ and $L \to \infty$ finishes the proof.
\end{proof}

\begin{lem}\label{lem: close ngbhs keta} Fix $x \in X \setminus \bigcup \limits_{p=1}^{k} \Per_p(T)$. Then for $\Leb$-a.e. $\alpha \in B_m(0,1)$ there exist $\eta = \eta(x, \alpha)$ and $\gamma = \gamma(x,\alpha)$ such that for every $y \in X$
\begin{equation}\label{eq: close ngbhs keta}
\|\phi_\alpha(y) - \phi_\alpha(x)\| \leq \gamma \Rightarrow y \in K_\eta(x).
\end{equation}
\end{lem}

\begin{proof}
Set $P = P_{k}(x) \setminus B(x, \xi(x)/2)$. By Lemma \ref{lem: P disjoint k}, for $\Leb$-a.e. $\alpha \in B_m(0,1)$ we have $\phi_\alpha(\{x\}) \cap \phi_\alpha(P) = \emptyset$. Fix such $\alpha$. As $P$ is compact and $\phi_\alpha$ is continuous, there exist $\gamma > 0$ and $0 < \eta < \xi(x)/2$ such that for every $y \in X$
\[ \|\phi_\alpha(x) - \phi_\alpha(y)\| \leq \gamma \Longrightarrow \dist(y, P) \geq \eta,\]
hence also
\[ \|\phi_\alpha(x) - \phi_\alpha(y)\| \leq \gamma \Longrightarrow y \in K_{\eta}(x).\]

Indeed, if $ \dist(y, P) \geq \eta$ and $\dist(y, P_{k}(x)) < \eta$, then there exists $z \in P_{k}(x) \setminus P \subset B(x, \xi(x)/2)$ such that $\|z - y\| < \eta < \xi(x)/2$, so $y \in B(x, \xi(x))$.
\end{proof}

Now we can complete the proof of Proposition \ref{prop: unpred}.

\begin{proof}[{\bf Proof of Proposition \ref{prop: unpred}}]
Fix $x \in X \setminus \bigcup \limits_{p=1}^{k} \Per_p(T)$. By Lemma \ref{lem: close ngbhs keta}, for $\Leb$-a.e. $\alpha \in B_m(0,1)$ there exist $\eta = \eta(x,\alpha) > 0$ and $n_2 = n_2(x, \alpha)$ such that for every $j \in \Z$ such that $|j| \geq n_2$ we have (applying \eqref{eq: close ngbhs keta} with $y = T^j x$, for $n_2$ such that $n_2^{-1/(k+\eps')} < \gamma(x,\alpha)$)
\[ \|\phi_\alpha(T^j x) - \phi_\alpha(x) \| \leq |j|^{-1/(k+\eps')} \Rightarrow T^j x \in K_\eta(x).\]
Therefore, by Lemma \ref{lem: keta unpred} we have that for $\Leb$-a.e. $\alpha$ there exists $n_0 = n_0(\alpha, x, \beta, L, \eps, \eps')$  such that \eqref{eq: unpred} holds for every $j \in \Z$ with $|j| \geq n_0$. As $\mu(\bigcup \limits_{p=1}^{k} \Per_p(T)) = 0$, Fubini's theorem finishes the proof of Proposition \ref{prop: unpred}.
\end{proof}

\subsection[Proof of Proposition \ref{prop: unrec}]{Proof of Proposition \ref{prop: unrec}}\label{subsec: proof unrec}

Now we shall continue with the proof of Proposition \ref{prop: unrec}. In what follows we assume that all the assumptions of Theorem \ref{thm: inverse FNN general} hold and that $T$ is a bijection. The super-polynomial decay of correlations will be used through the following lemma.

\begin{lem}\label{lem: spdc on balls}
	Fix $z \in X$. For every $j \in \Z$ and every  $r, \delta \in (0,1]$
	\[ \mu(B(z, r) \cap T^{-j}(B(z, \delta))) \leq \frac{4\theta(|j|)}{r\delta} + \mu(B(z, 2r))\mu(B(z, 2\delta)), \]
	where $\theta(\cdot)$ is as in the definition of super-polynomial decay of correlations (Definition \ref{defn: spdc}).
\end{lem}

\begin{proof}
	The proof is essentially the same as the proof of \cite[Lemma 6]{GalatoloHitting07} or \cite[Lemma 15]{RousseauSaussol10}. We include details for the convenience of the reader.
	
	Let $\eta_r : [0, \infty) \to \R$ be a $\frac{1}{r}$-Lipschitz function such that $\mathds{1}_{[0,r]} \leq \eta_r \leq \mathds{1}_{[0,2r]}$. Define $\psi_{z,r} : X \to \R$ as $\psi_{z,r}(x) = \eta_r(\|x - z\|)$. Note that $\|\psi_{z,r}\|_{\Lip} \leq  \|\psi_{z,r}\|_\infty + \Lip(\psi_{z,r}) \leq 1 + \frac{1}{r} \leq \frac{2}{r}$, as $r \leq 1$. By the super-polynomial decay of correlations we have
	\[\begin{split}
		\mu(B(z, r) \cap T^{-j}(B(z, \delta))) & \leq  \int \psi_{z,r}(x)\psi_{z, \delta}(T^j x)d\mu(x) \\
		& \leq \theta(|j|)\|\psi_{z,r}\|_{\Lip} \|\psi_{z,\delta}\|_{\Lip}  + \int \psi_{z,r}(x)d\mu(x) \int \psi_{z, \delta}(x)d\mu(x) \\
		& \leq \frac{4\theta(|j|)}{r\delta} + \mu(B(z, 2r))\mu(B(z, 2\delta)).
	\end{split} \]
\end{proof}

For $\rho, a,b,c > 0$ set
\begin{equation}\label{eq: G def}
G_{\rho,a,b,c} = \left\{ x \in X : \quad
\begin{aligned}
	&\underset{\stackrel{y \in X}{\|x - y\| \leq \rho}}{\forall} \sigma_k(D_{x,y}(k)) \geq \rho \|x - y\|, \\
	&\mu(B(x,r)) \leq r^{\hdim \mu - a}, \\
	&\mu(B(x,r)) \geq r^{N+b}, \\
	&\mu(B(x,4r)) \leq\mu(B(x,r)) r^{-c}
\end{aligned} \quad \text{ for every } 0 < r \leq \rho\right\}.
\end{equation}

\begin{lem}\label{lem: Grhoabc}
For every $a,b,c > 0$
\[ \lim \limits_{\rho \to 0} \mu(G_{\rho,a,b,c}) = 1. \]
\end{lem}

\begin{proof}
The proof follows from the following observations:
\begin{itemize}
	\item By \eqref{eq: sigma bound close}, for every $x \in X \setminus \bigcup \limits_{p=1}^{k} \Per_p(T)$ there exists $\rho = \rho(x)$ such that $\sigma_k(D_{x,y}(k)) \geq \rho  || x - y\|$ for every $y \in X$ such that $\|x-y\| \leq \rho$. Moreover, $\mu(\bigcup \limits_{p=1}^{k} \Per_p(T)) = 0$ (see the beginning of the proof of Theorem \ref{thm: inverse FNN general}).
	\item By Lemma \ref{lem: lower exact dim}, $\ld(\mu, x) = \hdim \mu$ for $\mu$-a.e.~$x$. Hence for every $a > 0$ and $\mu$-a.e. $x \in X$ there exists $\rho = \rho(x,a)$ such that $\mu(B(x,r)) \leq r^{\hdim \mu - a}$ for every $0 < r \leq \rho$.
	\item As $\mu$ is a measure on $\R^N$ we have $\ud(\mu, x) \leq N$ for $\mu$-a.e. $x$, and hence for every $b > 0$ and $\mu$-a.e. $x \in X$ there exists $\rho = \rho(x,b)$ such that $\mu(B(x,r)) \geq r^{N + b}$ for every $0 < r \leq \rho$.
	\item By \cite[Lemma 1]{barreira2001hausdorff}, $\mu$ is weakly diametrically regular (as a Borel probability measure in $\R^N$), hence for every $c > 0$ and $\mu$-a.e. $x \in X$ there exists $\rho = \rho(x,c)$ such that $\mu(B(x,4r)) \leq\mu(B(x,r)) r^{-c}$ for every $0 < r \leq \rho$.
\end{itemize}
\end{proof}

The main ingredient of the proof of Proposition \ref{prop: unrec} is the following inequality.

\begin{lem}\label{lem: main ineq spdc}
For every $\rho \in (0,1)$ and $p, a,b,c > 0$ such that $\hdim \mu > k + a$ and $p > k$ there exists $C > 0$ such that for every $0 < r \leq \delta < \rho/4$ and every $j \in \Z$
\[\begin{split}
\Leb & \otimes \mu  \left( \left\{ (\alpha, x) \in B_m(0,1) \times G_{\rho,a,b,c}: \|T^j x - x\| \leq \delta, \| \phi_\alpha(T^j x) - \phi_\alpha(x)\| \leq r\right\} \right) \\
& \leq C\left( r^{k-c}\delta^{d - a - k} + \theta(|j|)r^{k - N - b -1}\delta^{-k-1} + \theta(|j|)r^{k - N - b -1}\delta^{-p(1 + \frac{1}{k})}  + r^{k-N-b}\delta^{p(\frac{d - a}{k} - 1)} \right),
\end{split}\] 
where $d = \hdim \mu$.
\end{lem}

\begin{proof}
Denote $G = G_{\rho, a, b, c}$. Fix $0 < r \leq \delta < \rho/4$ and $j \in \Z$. By \cite[Lemma 2]{barreira2001hausdorff} there exists a set $E \subset G$ such that
\begin{equation}\label{eq: E prop}
\mu \left( G \setminus \bigcup \limits_{y \in E} B(y,r) \right) = 0 \text{ and the family } \{B(z, r/2) : z \in E\} \text{ is pairwise disjoint}. 
\end{equation}
Note that
\[ 1 \geq \sum \limits_{z \in E} \mu(B(z, r/2)) \geq 2^{-(N+b)}r^{N+b}\# E, \]
so
\begin{equation}\label{eq: E cardinality}
\# E \leq 2^{N+b}r^{-(N+b)}.
\end{equation}
Note the following: if $x \in G$ and $z \in E$ are such that $x \in B(z,r)$, then (recall that we assume $r \leq \delta$)
\[ \|T^j x - x\| \leq \delta \Rightarrow \|T^j x - z\| \leq 2\delta \]
and
\[ \| \phi_\alpha(T^j x) - \phi_\alpha(x)\| \leq r \Rightarrow \| \phi_\alpha(T^j x) - \phi_\alpha(z)\| \leq Mr, \]
where $M = 1 + \sup \limits_{\alpha \in B_m(0,1)} \Lip(\phi_\alpha)$.
Therefore
\[\begin{split}
	\Leb & \otimes \mu  \left( \left\{ (\alpha, x) \in B_m(0,1) \times G: \|T^j x - x\| \leq \delta, \| \phi_\alpha(T^j x) - \phi_\alpha(x)\| \leq r\right\} \right) \\
	& \leq \sum \limits_{z \in E} \Leb  \otimes \mu  \Big( \Big\{ (\alpha, x) \in B_m(0,1) \times B(z,r) : \\
	& \qquad \qquad \|T^j x - z\| \leq 2\delta, \| \phi_\alpha(T^j x) - \phi_\alpha(z)\| \leq Mr\Big\} \Big) \\
	& = \sum \limits_{z \in E} \int \limits_{B(z,r)} \mathds{1}_{B(z,2\delta)}(T^j x)\Leb \left( \{ \alpha \in B_m(0,1) : \|\phi_\alpha(T^j x) - \phi_\alpha(z) \| \leq M r \} \right) d\mu(x) = (*).
\end{split}\]
Applying Lemma \ref{lem: key_ineq_inter} and recalling that $2\delta < \rho$, we can bound further
\[\begin{split}
(*) &\leq C \sum \limits_{z \in E} \int \limits_{B(z,r)} \mathds{1}_{B(z,2\delta)}(T^j x)\frac{(Mr)^k }{\sigma_k(D_{T^j x,z}(k))^{k}} d\mu(x) \\
& \leq CM^k \rho^{-k} r^k  \sum \limits_{z \in E} \int \limits_{B(z,r)} \mathds{1}_{B(z,2\delta)}(T^j x)\|T^j x - z\|^{-k} d\mu(x) = (*).
\end{split}\]
Changing variables (see e.g.~\cite[Theorem 1.15]{mattila}) one gets further
\[\begin{split}
(*) & =  CM^k \rho^{-k} r^k  \sum \limits_{z \in E} \int \limits_{0}^\infty \mu \left( \left\{ x \in B(z,r) \cap T^{-j}(B(z,2\delta)) : \|T^j x - z\|^{-k} \geq t  \right\} \right) dt \\
& = CM^k \rho^{-k} r^k  \sum \limits_{z \in E} \int \limits_{0}^\infty \mu \left( \left\{ x \in B(z,r) \cap T^{-j}(B(z,2\delta)) : \|T^j x - z\| \leq t^{-1/k}  \right\} \right) dt = (*).
\end{split}\]
We split the last expression into three terms:
\[\begin{split}
(*)  = &\  CM^k \rho^{-k} r^k  \sum \limits_{z \in E} \int \limits_{0}^{(2\delta)^{-k}} \mu \left( B(z,r) \cap T^{-j}(B(z,2\delta)) \right) dt \\
& + CM^k \rho^{-k} r^k  \sum \limits_{z \in E} \int \limits_{(2\delta)^{-k}}^{(2\delta)^{-p}} \mu \left( B(z,r) \cap T^{-j}(B(z,t^{-1/k})) \right) dt \\
& + CM^k \rho^{-k} r^k  \sum \limits_{z \in E} \int \limits_{(2\delta)^{-p}}^{\infty} \mu \left( B(z,r) \cap T^{-j}(B(z,t^{-1/k})) \right) dt \\
& =: I + II + III.
\end{split}\]
For the first one we employ Lemma \ref{lem: spdc on balls} (note that the integrand does not depend on $t$) and use \eqref{eq: E prop}, \eqref{eq: E cardinality} together with the definition of $G = G_{\rho,a,b,c}$ and $T$-invariance of $\mu$ to obtain
\[\begin{split}
I & \leq CM^k \rho^{-k} r^k \sum \limits_{z \in E} (2\delta)^{-k} \left( \frac{4\theta(|j|)}{2r\delta} + \mu ( B(z,2r)) \mu (T^{-j}(B(z,4\delta))) \right) \\
& \overset{\eqref{eq: G def}}{\leq} CM^k2^{-k} \rho^{-k} r^k \delta^{-k} \left( \# E \frac{4\theta(|j|)}{2r\delta} + 2^cr^{-c}\sum \limits_{z \in E} \mu ( B(z,r/2)) \mu (B(z,4\delta)) \right) \\
& \overset{\eqref{eq: G def},\, \eqref{eq: E cardinality}}{\leq} CM^k2^{N+b-k} \rho^{-k} r^k \delta^{-k} \left( 2\theta(|j|)r^{-(N+b + 1)}\delta^{-1} + 2^cr^{-c}\sum \limits_{z \in E} \mu ( B(z,r/2))  (4\delta)^{d-a}\right) \\
& \overset{\eqref{eq: E prop}}{\leq}  CM^k2^{N+b-k} \rho^{-k} r^k \delta^{-k} \left( 2\theta(|j|)r^{-(N+b + 1)}\delta^{-1} + 2^cr^{-c} (4\delta)^{d-a}\right) \\
& = C_{\rho, a, b, c} \left(r^{k-c} \delta^{d-a-k} + \theta(|j|)r^{k - N - b - 1}\delta^{-k -1} \right).
\end{split}\]
The second term is bounded similarly:
\[
\begin{split}
II & \leq CM^k \rho^{-k} r^k  \sum \limits_{z \in E} \int \limits_{(2\delta)^{-k}}^{(2\delta)^{-p}} \left( \frac{4\theta(|j|)}{rt^{-{1/k}}} + \mu(B(z, 2r))\mu(B(z, 2t^{-1/k})) \right) dt \\
& \overset{\eqref{eq: G def},\, \eqref{eq: E cardinality}}{\leq} CM^k \rho^{-k} r^k \bigg( 2^{N+b+2} \theta(|j|)r^{-N-b-1} \int \limits_{0}^{(2\delta)^{-p}} t^{1/k}dt \\
& \qquad \qquad + 2^cr^{-c}\sum \limits_{z \in E} \mu ( B(z,r/2)) \int \limits_{(2\delta)^{-k}}^{(2\delta)^{-p}} \mu(B(z, 2t^{-1/k}))dt\bigg) \\
& \overset{\eqref{eq: G def},\, \eqref{eq: E prop}}{\leq}  CM^k \rho^{-k} r^k \left( 2^{N+b+2}  \theta(|j|)r^{-N-b-1} \frac{(2\delta)^{-p(1 + \frac{1}{k})}}{1+ \frac{1}{k}} + 2^cr^{-c}\int \limits_{(2\delta)^{-k}}^{\infty} 2^{d-a}t^{-(d-a)/k}dt\right) \\
& =  CM^k \rho^{-k} r^k \left( \frac{2^{N+b+2-p(1+\frac{1}{k})}}{1+ \frac{1}{k}} \theta(|j|)r^{-N-b-1} \delta^{-p(1 + \frac{1}{k})} +  2^{d-a}2^cr^{-c}\frac{(2\delta)^{ d - a - k}}{\frac{d-a}{k} - 1}\right) \\
& \leq C_{\rho,a,b,c,p} \left(r^{k-c}\delta^{d-a-k} + \theta(|j|)r^{k - N - b - 1}\delta^{-p(1 + \frac{1}{k})} \right).
\end{split}
\]
Finally, for the third term we have by \eqref{eq: E cardinality}
\[
\begin{split}
III &\leq  CM^k \rho^{-k} r^k  \sum \limits_{z \in E} \int \limits_{(2\delta)^{-p}}^{\infty} \mu \left(B(z,t^{-1/k}) \right) dt  \leq  C 2^{N+b}M^k \rho^{-k} r^{k-N-b}  \int \limits_{(2\delta)^{-p}}^{\infty} t^{-\frac{d-a}{k}} dt \\
& = C 2^{N+b} M^k \rho^{-k} r^{k-N-b} \frac{(2\delta)^{p(\frac{d-a}{k} - 1)}}{\frac{d-a}{k} - 1}  = C_{\rho,a,b,c,p} r^{k - N - b}\delta^{p(\frac{d-a}{k} - 1)}.
\end{split}
\]
Combining the three bounds finishes the proof.
\end{proof}

Now we can prove Proposition \ref{prop: unrec} (which, as already explained, combined with Proposition \ref{prop: unpred} yields Theorem \ref{thm: inverse FNN general}).

\begin{proof}[{\bf Proof of Proposition \ref{prop: unrec}}]
Fix $\eps>0$. We shall prove that for every small enough $\eps'>0$, there exist $a,b,c>0$ such that for arbitrarily small $\rho$, setting $G = G_{\rho,a,b,c}$, we have that the series
\begin{equation}\label{eq: sum to converge}
\begin{split}
\sum \limits_{j \in \Z} \Leb \otimes \mu \bigg( \bigg\{ (\alpha, x) \in B_m(0,1) \times G : \|&\phi_\alpha(T^j x) - \phi_\alpha(x) \| \leq |j|^{-\frac{1}{k + \eps'}} \\
& \text{ and }\|T^j x - x \| \leq |j|^{-\eps}\bigg\} \bigg)
\end{split}
\end{equation}
converges. This will finish the proof by applying the Borel--Cantelli lemma, Lemma \ref{lem: Grhoabc} and Fubini's theorem.

We can assume w.l.o.g. that $\eps < 1/k$. Given $\eps'>0$ (to be determined later) and $\rho$, let $M \in \N$ be such that
\[ |j| \geq M \Rightarrow |j|^{-\frac{1}{k+\eps'}} < \rho / 4, |j|^{-\eps} < \rho/4.\]
As we assume $\eps < 1/k$, we have (for $\eps'$ small enough depending on $\eps$) that $|j|^{- \frac{1}{k+\eps'}} \leq |j|^{-\eps}$ for all $j \in \Z \setminus \{0\}$, hence we can apply Lemma \ref{lem: main ineq spdc} with $\delta = |j|^{-\eps}$ and $r = |j|^{-1/(k+\eps')}$ to obtain for every $a,b,c > 0$ and $p > k$, with $d = \hdim \mu > k+a$
\[
\begin{split}
\sum  \limits_{j \in \Z} &  \Leb \otimes \mu  \left( \left\{ (\alpha, x) \in B_m(0,1) \times G : \|\phi_\alpha(T^j x) - \phi_\alpha(x) \| \leq |j|^{-\frac{1}{k + \eps'}},\ \|T^j x - x \| \leq |j|^{-\eps}\right\} \right) \\
& \leq 2M-1 +  \sum \limits_{|j| \geq M} C_{\rho,a,b,c,p} \Big( |j|^{-\frac{k-c}{k+\eps'}}|j|^{-\eps(d-a-k)} + \theta(|j|)|j|^{-\frac{k-N-b-1}{k+\eps'}}|j|^{\eps(k+1)}  \\
& \qquad  \qquad \qquad\qquad\qquad\qquad + \theta(|j|)|j|^{-\frac{k-N-b-1}{k+\eps'}}|j|^{\eps p (1 + \frac{1}{k})}  + |j|^{-\frac{k-N-b}{k+\eps'}}|j|^{-\eps p(\frac{d-a}{k} - 1)} \Big).
\end{split}\]
Note the following:
\begin{itemize}
\item $\sum \limits_{|j| \geq M} \left( \theta(|j|)|j|^{-\frac{k-N-b-1}{k+\eps'}}|j|^{\eps(k+1)} +  \theta(|j|)|j|^{-\frac{k-N-b-1}{k+\eps'}}|j|^{\eps p (1 + \frac{1}{k})} \right) < \infty$ for any choice of parameters $\eps, \eps', a,b,c,p$ as $\theta(|j|)$ decays super-polynomially in $|j|$.
\item $\sum \limits_{|j| \geq M} |j|^{-\frac{k-c}{k+\eps'}}|j|^{-\eps(d-a-k)} < \infty$ if $\eps(d-a-k) > 1 - \frac{k-c}{k+\eps'}$. For fixed $\eps>0$ and $a > 0$ this can be obtained by choosing $c$ and $\eps'$ small enough.
\item $\sum \limits_{|j| \geq M} |j|^{-\frac{k-N-b}{k+\eps'}}|j|^{-\eps p(\frac{d-a}{k} - 1)} < \infty$ if $\eps p \left(\frac{d-a}{k} - 1\right) > 1 + \frac{N+b - k}{k + \eps'}$. For fixed $\eps, \eps',a,b,c$ this can be obtained by choosing $p$ large enough.
\end{itemize}
Therefore, for given $\eps < \frac{1}{k}$, for every $\eps' > 0$ small enough one can choose $a,b,c,p$ in such a way that the series \eqref{eq: sum to converge} converges for every $\rho \in (0,1)$. This finishes the proof.
\end{proof}

\section[Applications to H\'enon and Lorenz attractors]{Applications to H\'enon and Lorenz attractors}\label{sec: applications}

{The H\'enon map and the Lorenz system are standard test cases for the FNN algorithm. For noise-free time series, the proportion of false nearest neighbours is reported to drop to (almost) zero at $k = 2$ for the H\'enon map and at $k = 3$ for the Lorenz system, see \cite{TISEAN99} and \cite[Section~4.2.1]{Abarbanel_book}. In this section we apply Theorems \ref{thm: FNN diffeo} and \ref{thm: inverse FNN diffeo} to show that for typical observables and trajectories, the FNN algorithm with the stopping criterion (4') and a large threshold $R_{\tol}$ returns
\begin{enumerate}[(i)]
\item\label{it: henon upper} $k \leq 2$ for the H\'enon map for all dissipative parameters and every invariant measure on the attractor,
\item\label{it: henon exact} $k=2$ for the Benedicks--Carleson parameters of the H\'enon map, with respect to the natural (SRB) measure on the attractor and for Lebesgue-almost every initial point in its basin of attraction,
\item\label{it: lorenz upper} $k \leq 3$ for the Lorenz system, for all (dissipative) parameters and every invariant measure.
\end{enumerate}

Both systems are defined on Euclidean spaces rather than on compact manifolds. As all the orbits considered below stay in a compact region, this is a technical matter: after modifying the systems away from this region, one obtains diffeomorphisms of the spheres $S^2$ and $S^3$, to which the results of the previous sections apply. We shall not go into details.

The assumption $\mH^{\dim M}(\supp \mu) = 0$ in \ref{it: A2} holds for dissipative systems. Indeed, let $\mu$ be a compactly supported $T$-invariant Borel probability measure and $X = \supp \mu$. Then $X$ is compact, $\mu(X) = 1$ and $T(X) = X$. If $|\det DT| \leq b < 1$ on $X$, then it follows directly that $X$ has zero Lebesgue measure, hence $\mH^{\dim M}(\supp \mu) = 0$.

\subsection{H\'enon maps}\label{subsec: henon}

For $a, b \in \R$ with $b \neq 0$, the H\'enon map \cite{Henon76} is
\[ T_{a,b} : \R^2 \to \R^2, \qquad T_{a,b}(x,y) = (1 - a x^2 + y,\ b x). \]
It is a polynomial diffeomorphism of $\R^2$ with constant Jacobian $\det DT_{a,b} \equiv -b$. Hence for $0 < |b| < 1$, we have $\mH^2(\supp \mu) = 0$ for every compactly supported $T_{a,b}$-invariant Borel probability measure $\mu$, so \ref{it: A2} holds with $k = 2$ and Theorem \ref{thm: FNN diffeo} can be applied, proving \ref{it: henon upper} for $\mu$-typical points.
For a lower bound, we use Theorem \ref{thm: inverse FNN diffeo}, which requires super-polynomial decay of correlations and $\hdim \mu > 1$. Both are known for the Benedicks--Carleson parameters. The results that we invoke below are summarized in the survey \cite{Benedicks02}. By \cite{BC91}, there exists $b_0 > 0$ such that for every $0 < b < b_0$ there exists a certain set $\Delta_b$ of positive one-dimensional Lebesgue measure of parameters $a < 2$. For $a \in \Delta_b$ the map $T = T_{a,b}$ has a compact attractor $X$ with $T(X) = X$ on which the dynamics is chaotic (has sensitive dependence on initial conditions). The set
\[ \mathcal{A} = \{ (a,b) : 0 < b < b_0,\ a \in \Delta_b \} \]
has positive two-dimensional Lebesgue measure, and its elements are called the \emph{Benedicks--Carleson parameters}. By \cite{BY93}, for $(a,b) \in \mathcal{A}$ the attractor $X$ supports a unique SRB measure $\mu$, with Lyapunov exponents $\lambda_1 > 0 > \lambda_2$, and $(T, \mu)$ has exponential decay of correlations for H\"older observables \cite{BenedicksYoung00}, so also a super-polynomial decay of correlations in the sense of Definition \ref{defn: spdc}. The measure $\mu$ is the natural measure on $X$: $\frac{1}{n}\sum_{i=0}^{n-1} \delta_{T^i x} \to \mu$ weakly-$*$ for Lebesgue-almost every point $x$ of the basin of attraction of $X$ \cite{BV01}, so it describes the statistics of a physically typical initial condition. Moreover, by \cite[Section~4.5]{BY93},
\[ \hdim \mu = 1 + \frac{\lambda_1}{\lambda_1 + |\log b|} \in (1,2). \]
We can therefore apply Theorem \ref{thm: inverse FNN diffeo} with $k = 1$, proving \ref{it: henon exact} for $\mu$-typical points.

Let us now briefly describe how one can also use our results to cover Lebesgue-typical points in the basin of attraction. We provide only a sketch and leave the details to the interested reader. For $(a,b)\in\mathcal A$, Lebesgue-almost every point $x$ in the basin of $X$ lies on the stable manifold of a $\mu$-typical point $z\in X$, so that $\|T^i x-T^i z\|=O(e^{-ci})$ for some $c>0$. This is the basin theorem of Benedicks and Viana \cite{BV01} (see also \cite[Proposition~4.2 and its proof]{WangYoung13}). It follows from the proof of Theorem~\ref{thm: inverse FNN diffeo} at delay length one that the lower bound $k\geq 2$ carries over from a typical orbit of $z$ to any orbit tracking it exponentially. We may discard finitely many initial iterates without changing the limiting proportions. For a proportion of indices tending to one, exponential tracking preserves the polynomial bounds used in applying Propositions~\ref{prop: unpred} and~\ref{prop: unrec}, namely the upper bound on the observed nearest-neighbour distance and the lower bound on the observed distance between the same two orbit points one iterate later. The same tracking argument gives the upper bound $k\leq 2$. Indeed, for a prevalent observable $h$, the measure $\phi_{h,2}\mu$ is exact dimensional with $\hdim(\phi_{h,2}\mu)=\hdim\mu>0$ by \cite[Corollary~A.2]{BGSLimits}, so the theorem of Rousseau and Saussol \cite[Theorem~5]{RousseauSaussol10} gives the precise polynomial recurrence rate for the two-delay observations on the attractor.

Note that the classical parameters $a = 1.4$, $b = 0.3$ considered by H\'enon \cite{Henon76} are not known to be Benedicks--Carleson parameters. For them, the existence of an SRB measure is not known, and periodic sinks were found for parameters extremely close to $(1.4, 0.3)$ \cite{GaliasTucker15}. The upper bound $k \leq 2$ applies to them, as to all parameters with $0 < |b| < 1$.

\subsection{The Lorenz system}\label{subsec: lorenz}

For $\sigma, \beta, \rho > 0$, the Lorenz system \cite{Lorenz63} is given by the 3-dimensional ODE
\begin{equation}\label{eq: lorenz} \dot x = \sigma(y - x), \qquad \dot y = x(\rho - z) - y, \qquad \dot z = xy - \beta z. \end{equation}
Let $(\varphi_t)_{t \in \R}$ be its flow on $\R^3$. There exists a compact set $B \subset \R^3$ such that $\varphi_t(B) \subset B$ for $t \geq 0$ and every trajectory eventually enters $B$ \cite{Lorenz63}. As the vector field of the equation \eqref{eq: lorenz} has negative divergence $-(\sigma + 1 + \beta)$, the system is dissipative. Therefore, for $T = \varphi_\tau$ being the time-$\tau$ map for some $\tau > 0$, we have that every $T$-invariant Borel probability measure is compactly supported inside $B$ and $\det D\varphi_\tau < 1$ gives $\mH^3(\supp \mu) = 0$. Therefore \ref{it: A2} holds with $k = 3$ and Theorem \ref{thm: FNN diffeo} can be applied for all $\sigma, \beta, \rho, \tau > 0$ and every $\varphi_\tau$-invariant Borel probability measure $\mu$, proving \ref{it: lorenz upper}.

To obtain the matching lower bound $k = 3$ for the classical parameters $\sigma = 10$, $\beta = 8/3$, $\rho = 28$ from Theorem \ref{thm: inverse FNN diffeo}, one needs super-polynomial decay of correlations and $\hdim \mu > 2$.  The existence of the SRB measure $\mu$ on the Lorenz attractor  was established first by Tucker \cite{Tucker02} (see also \cite{APPV09}) and its exponential decay of correlations for H\"older observables was obtained in \cite{AM16, AM17}. Numerically, the dimension of the Lorenz attractor is estimated as $\approx 2.06$ \cite{Viswanath04}, but the bound $\hdim \mu > 2$ does not seem to have been established. Rigorously, geometric Lorenz attractors (a simplified model of the Lorenz attractor) contain horseshoes of Hausdorff dimension larger than $2$ \cite{MPR20, RegoVarandas26}, but this concerns sets rather than the measure $\mu$.%
}

{\section*{Acknowledgements}

This work originated in a long-standing joint project with Krzysztof Bara\'nski and Yonatan Gutman on the probabilistic Takens theory. The author is very grateful to them for many inspiring discussions over the years.

This work was partially supported by the Simons Foundation grant (award no.~SFI-MPS-T-Institutes-00010825) and from State Treasury funds as part of a task commissioned by the Minister of Science and Higher Education under the project ``Organization of the Simons Semesters at the Banach Center -- New Energies in 2026--2028'' (agreement no.~MNiSW/2025/DAP/491).

\vspace{0.5cm}
\begin{center}
	\includegraphics[width=0.8\textwidth]{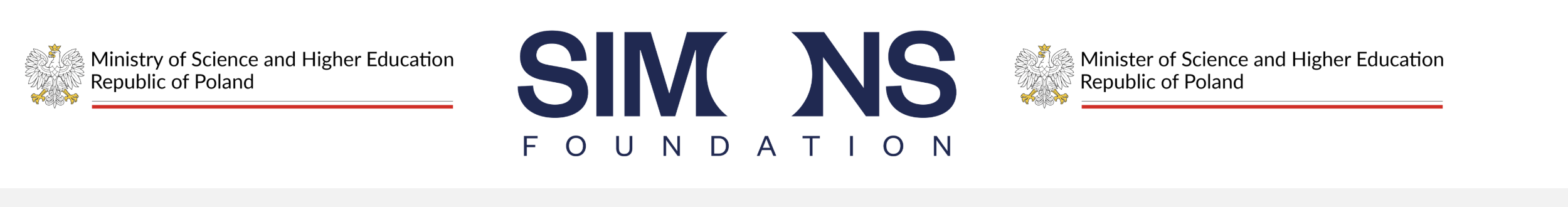}
\end{center}

\section*{Declaration on the use of AI tools}

During the preparation of this work the author used an AI model for typesetting, language editing, literature search and editing of the text based on the author's instructions. In particular, some of the references in Sections 1 and 5 were proposed by the AI model and have been verified by the author. All proofs and mathematical ideas were developed by the author without the use of any AI tools. The author takes full responsibility for the content of the paper.
}

\bibliographystyle{alphaurl}
\bibliography{FNN}

\end{document}